
\input epsf.tex
\input psfrag
\newdimen\epsfxsize
\newdimen\epsfysize
\magnification=\magstephalf 
\baselineskip=18pt




\font\tenmsb=msbm10
\font\sevenmsb=msbm7
\font\fivemsb=msbm5
\newfam\msbfam
\textfont\msbfam=\tenmsb
\scriptfont\msbfam=\sevenmsb
\scriptscriptfont\msbfam=\fivemsb
\def\Bbb#1{\fam\msbfam#1}

\font\tenmsa=msam10
\font\sevenmsa=msam7
\font\fivemsa=msam5
\newfam\msafam
\textfont\msafam=\tenmsa
\scriptfont\msafam=\sevenmsa
\scriptscriptfont\msafam=\fivemsa
\mathchardef\blacksquare="0\number\msafam04

\def\C{{\Bbb C}}
\def\D{{\Bbb D}}
\def\H{{\Bbb H}}
\def\R{{\Bbb R}}





\def\diam{{\mathop {\rm diam}}}
\def\dist{{\mathop {\rm dist}}}

\def\Im{{\mathop {\rm Im}}}
\def\Re{{\mathop {\rm Re}}}


\def\proof{\bigskip\noindent {\bf Proof.} }
\def\endproof{\hfill $\blacksquare$\bigskip}
\def\endpf{{\qquad {\blacksquare}}}

\def\f{{\varphi}}
\def\newchi{\raise2pt\hbox{$\chi$}}
\def\newdag{\raise4pt\hbox{\dag}}
\def\eps{\varepsilon}

\def\Om{{\Omega}}
\def\om{{\omega}}
\def\b{{\partial }}
\def\bD{{\partial \D}}




%
\centerline{Convergence of the Zipper algorithm for conformal mapping}
\bigskip
\centerline{Donald E. Marshall\newdag and Steffen Rohde\newdag}
\centerline{Department of Mathematics}
\centerline{University of Washington}
\footnote{}{\newdag The authors are supported in part by NSF grants DMS-0201435
and DMS-0244408.} 

\bigskip
\bigskip
\centerline{\bf Abstract}

In the early 1980's an elementary algorithm for computing conformal
maps was discovered by R. K\"uhnau and the first author. The algorithm is fast
and accurate, but convergence was not known. Given points $z_0,\dots,z_n$
in the plane, 
the algorithm computes an explicit conformal map of the unit disk onto
a region bounded by a Jordan curve $\gamma$ with $z_0,\dots,z_n \in
\gamma$.
We prove convergence  for Jordan regions 
in the sense of uniformly close
boundaries, and give corresponding uniform estimates on the closed
region and the closed disc for the mapping functions and their
inverses.  Improved estimates are obtained if
the data points lie on a $C^1$ curve or a $K-$quasicircle.
The algorithm was discovered as an
approximate method for conformal welding, however it can also be
viewed as a discretization of the L\"owner differential equation.

\bigskip

\S {\bf 0. Introduction}
\bigskip

Conformal maps have useful applications to problems in physics, engineering
and mathematics, but how do you find a conformal map say of the
upper half plane $\H$ to a complicated region? Rather few maps can be
given explicitly by hand, so that a computer must be used to find
the map approximately. One
reasonable way to describe a region numerically is to give a large
number of points on the boundary. One way to say that a computed
map defined on $\H$ is ``close'' to a map to the region is to require that the
boundary of the image be uniformly close to the polygonal curve
through the data points.  Indeed, the only information we may have
about the boundary of a region are these data points.

\bigskip
\psfrag{vp}{$\varphi$}
\psfrag{vpi}{$\varphi^{-1}$}
\psfrag{z0}{$z_0$}
\psfrag{z1}{$z_1$}
\psfrag{z2}{$z_2$}
\psfrag{zn}{$z_n$}
\psfrag{Om}{$\Omega$}
\psfrag{H}{$\H$}
\centerline{\includegraphics[height=2.25in]{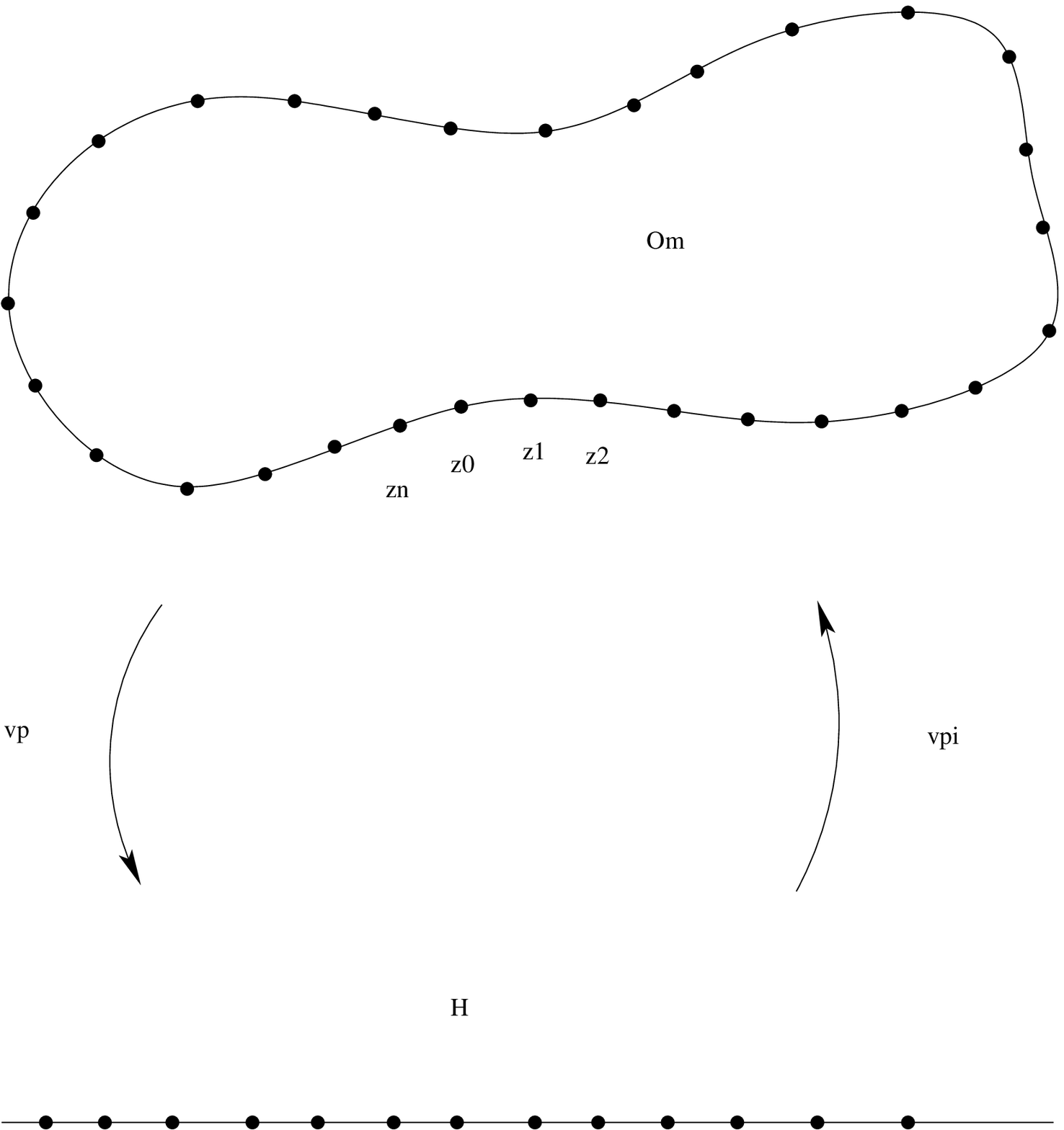}}
\nobreak
\centerline{{\bf Figure 1.}}
\vskip 0.5truein

In the early 1980's an elementary algorithm was
discovered independently by R. K\"uhnau [K] and the first author. The
algorithm is fast and accurate, but convergence was not known. 
The purpose of this paper is to prove convergence in the sense of
uniformly close boundaries, and discuss related
numerical issues. 
One important aspect of the algorithm that sets it apart from others:
in many applications both the conformal map and its inverse are
required; this algorithm finds both simultaneously.

The algorithm
can be viewed as a discretization of the Loewner differential equation,
or as an approximate solution to a conformal welding problem.
The approximation to the conformal map is obtained as a composition
of conformal maps onto slit halfplanes. Depending on the type of slit
(hyperbolic geodesic, straight line segment or circular arc) we actually 
obtain different versions of this algorithm. These are described in Section 1.

We then focus our attention on the ``geodesic algorithm'' and study its behaviour in
different situations. The easiest case is discussed in Section 2: If the data points
$z_0,z_1,...$ are the consecutive contact points of a chain of disjoint discs (see Figures 7 and 8 below), then a simple but very useful reinterpretation of the algorithm, together
with the hyperbolic convexity of discs in simply connected domains (J\o rgensen's theorem),
implies that the curve produced by the algorithm is
confined to the chain of discs (Theorem 2.2). 
One consequence is that for any bounded simply connected domain $\Om$,
the geodesic algorithm can be used to compute  a conformal map 
to a Jordan region $\Om_c$ (``c'' for computed) 
so that the Hausdorff distance between 
$\b \Om$ and $\b\Om_c$ is as small as desired (Theorem 2.4).

In Section 3, we describe an extension of the ideas of Section 2 that applies to a variety
of domains such as smooth domains or quasiconformal discs with small constants, with better 
estimates. For instance, if $\b \Om$ is a $C^1$ curve, then the geodesic
algorithm can be used to compute a conformal map to a Jordan region $\Om_c$ 
with $\b\Om_c\in C^1$ so that the boundaries are uniformly close and 
so that the 
unit tangent vectors are uniformly close (Theorem 3.10).
The heart of the convergence proof in these cases is the
technical ``self-improvement'' 
Lemmas 3.5 and 3.6. In fact, this approach constituted our first convergence proof.

The basic conformal maps and their inverses used in the geodesic
algorithm are given in terms of linear fractional transformations,
squares and square roots. The slit and zipper algorithms use
elementary maps whose inverses cannot be written in terms of
elementary maps. Newton's method, however, converges so rapidly that
it provides virtually a formula for the inverses. In Section 4 we
discuss how to apply variants of Newton's method by dividing the plane 
into four
regions, and prove quadratic convergence in one region.  We plan to 
address the convergence of the slit and zipper variants of the
algorithm in a forthcoming paper.

In Sections 5 and 6, we show how estimates on the distance between
boundaries of Jordan regions gives estimates on the uniform distance between
the corresponding conformal maps to $\D$,
and 
apply these estimates to obtain bounds for the convergence of the conformal
maps produced by the algorithm. We summarize some of our results as follows:
If $\partial \Omega$ is contained in a chain of discs of radius $\leq \epsilon$ with the data
points being the contact points of the discs, or
if $\partial \Omega$ is a $K$-quasicircle with $K$ close to one and if the data points are consecutive
points on $\partial \Omega$ of distance comparable to $\epsilon$,
then the Hausdorff distance between $\b\Om$ and the boundary of the
domain computed by the geodesic algorithm, $\b\Om_c$, is at most
$\eps$
and the conformal maps $\f,\f_c$ onto $\D$ satisfy
$$\sup_{\Omega\cap \Omega_c} |\f - \f_c| \leq C \epsilon^p,$$
where any $p<1/2$ works in the disc-chain case, and $p$ is close to 1 if $K$ is close to one.
In the case of quasicircles, we also have
$$\sup_{\D} |\f^{-1} - \f_c^{-1}| \leq C \epsilon^p$$
with $p$ close to one. Better estimates are obtained for regions
bounded by smoother Jordan curves.

Section 7 contains a brief discussion of numerical results.
The Appendix has a simple self-contained proof of J\o rgensen's
theorem.

The first author would like to express his deep gratitude to L. Carleson
for our exciting investigations at Mittag-Leffler Institute 1882-83
which led to the discovery of the zipper algorithms.

\bigskip
\bigskip
\S {\bf 1. Conformal mapping algorithms}
\bigskip

\noindent {\bf The Geodesic Algorithm}
\bigskip
The most elementary version of the conformal mapping algorithm is
based on the simple map $f_a: \H\setminus \gamma\longrightarrow \H$ where
$\gamma$ is an arc of a circle from $0$ to $a\in \H$ 
which is orthogonal to $\R$ at $0$. This map can
be realized by a composition of a linear fractional transformation,
the square and
the square root map as illustrated in Figure~2. The 
orthogonal circle also meets $\R$ orthogonally at a point 
$b=|a|^2/\Re a$ and 
is illustrated by a dashed curve in Figure~2.

\vskip 0.5truein
\psfrag{Hmg}{$\H\setminus\gamma$}
\psfrag{H}{$\H$}
\psfrag{g}{$\gamma$}
\psfrag{a}{$a$}
\psfrag{b}{$b$}
\psfrag{c}{$c$}
\psfrag{mc}{$-c$}
\psfrag{ic}{$i c$}
\psfrag{c2}{$c^2$}
\psfrag{0}{$0$}
\psfrag{fa}{$f_a$}
\psfrag{l1}{$\displaystyle{{z}\over{1-z/b}}$}
\psfrag{sq}{$z^2+c^2$}
\psfrag{sqrt}{$\sqrt{z}$}
\centerline{\includegraphics[height=2.75in]{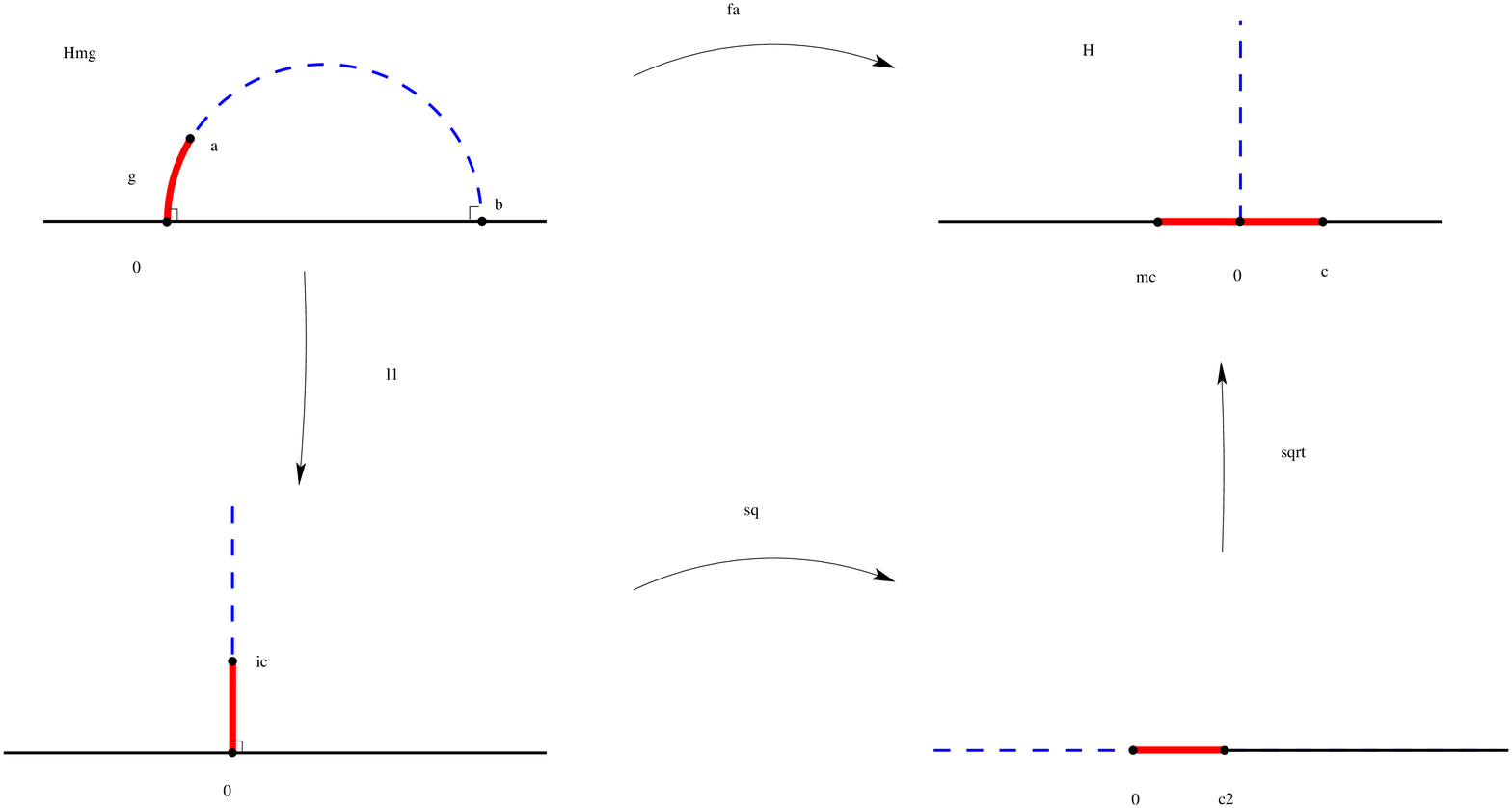}}
\nobreak
\centerline{{\bf Figure 2.} The basic map $f_a$.}
\vskip 0.5truein

\noindent In Figure 2, ~$c=|a|^2/\Im a$.
Observe that the arc $\gamma$ is opened to two adjacent intervals at
$0$ with $a$, the tip of $\gamma$, mapped to $0$.  
The inverse $f_a^{-1}$ can be easily found by composing
the inverses of these elementary maps in the reverse order.

Now suppose that $z_0, z_1, \dots, z_n$ are points in the plane. The
basic maps $f_a$ can be used to compute a conformal map of $\H$ onto a
region $\Om_c$ bounded by a Jordan curve which passes through the data points
as illustrated in Figure~3. 

\vskip 0.5truein
\psfrag{z0}{$z_0$}
\psfrag{z1}{$z_1$}
\psfrag{z2}{$z_2$}
\psfrag{z3}{$z_3$}
\psfrag{zn}{$z_n$}
\psfrag{ze2}{$\zeta_2$}
\psfrag{ze3}{$\zeta_3$}
\psfrag{ze4}{$\zeta_4$}
\psfrag{p1}{$\varphi_1=\displaystyle{i\sqrt{{(z-z_1)}/{(z-z_0)}}}$}
\psfrag{p2}{$\varphi_2=f_{\zeta_2}$}
\psfrag{p3}{$\varphi_3=f_{\zeta_3}$}
\psfrag{pn}{$\varphi_n=f_{\zeta_n}$}
\psfrag{pnp}{$\varphi_{n+1}=\displaystyle{-\biggl({z\over{1-z/\zeta_{n+1}}}\biggr)^2}$}
\psfrag{H}{$\H$}
\psfrag{0}{$0$}
\psfrag{Om}{$\Omega_c$}
\psfrag{zep}{$\zeta_{n+1}$}
\centerline{\includegraphics[height=4.0in]{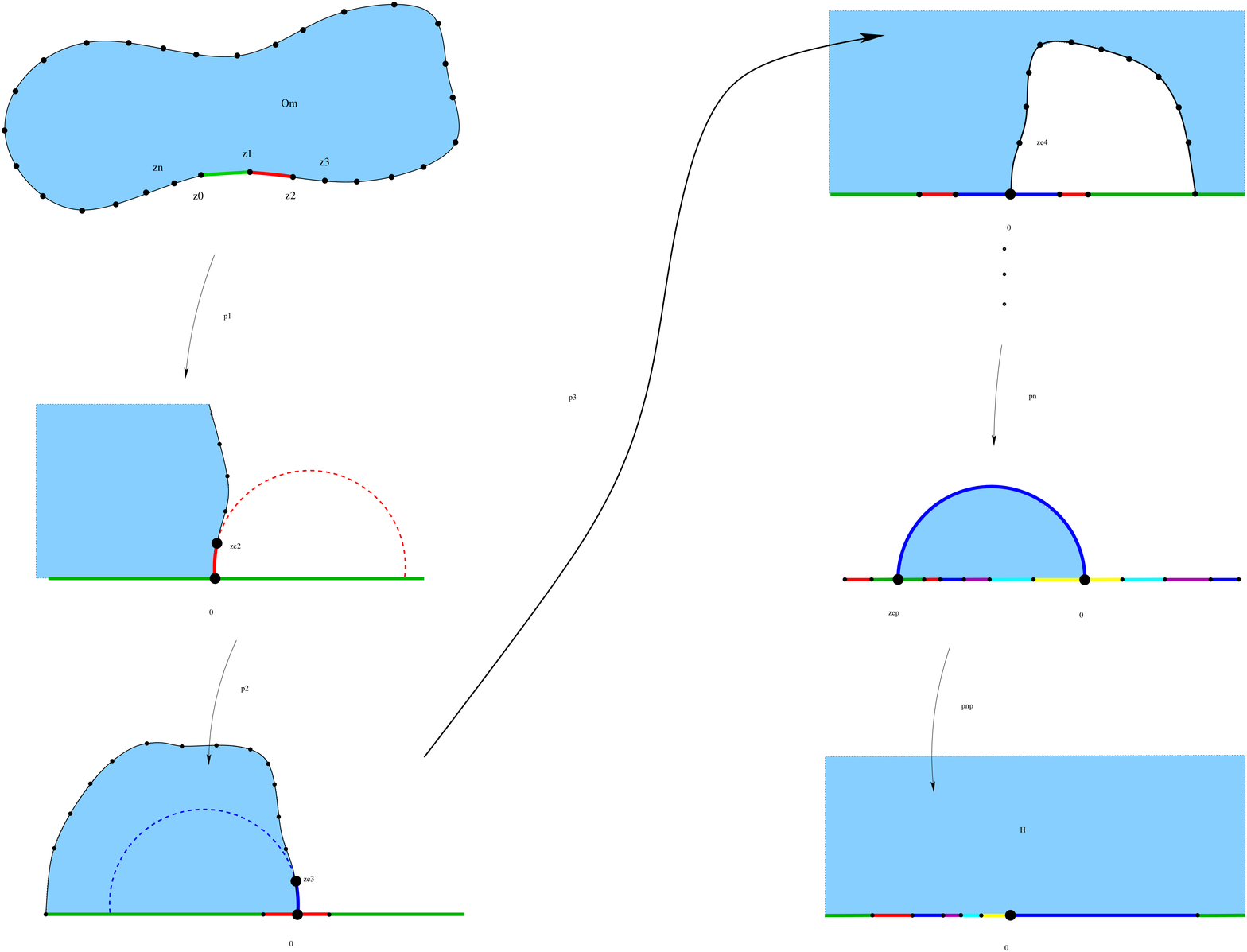}}
\nobreak
\centerline{{\bf Figure 3.} The Geodesic Algorithm.}
\vskip 0.5truein

The complement in the extended plane of the line segment from $z_0$ to
$z_1$ can be mapped onto $\H$ with the map 
$$\varphi_1(z)= i \sqrt{{z-z_1}\over{z-z_0}}$$
and $\varphi_1(z_1)=0$ and $\varphi_1(z_0)=\infty$. Set
$\zeta_2=\varphi_1(z_2)$ and $\varphi_2= f_{\zeta_2}$. 
Repeating this process, define
$$\zeta_k=\varphi_{k-1}\circ\varphi_{k-2}\circ\dots\circ\varphi_1(z_{k})$$
and
$$\varphi_k=f_{\zeta_k}.$$
for $k=2, \dots, n$. Finally, map a half-disc to $\H$ by letting 
$$\zeta_{n+1}=\varphi_n\circ\dots\circ\varphi_1(z_0)\in \R$$
be the image of $z_0$ and set
$$\varphi_{n+1}=
\pm\biggl({{z}\over{1-{z/{\zeta_{n+1}}}}}\biggr)^2$$
The $+$ sign is chosen in the definition of $\varphi_{n+1}$ if the data points
have negative winding number (clockwise) around an interior point of
$\b\Om$, and otherwise the $-$ sign is chosen.
Set
$$\varphi=\varphi_{n+1}\circ\varphi_{n}\circ\dots\circ\varphi_2\circ\varphi_1$$
and
$$\varphi^{-1}=\varphi_1^{-1}\circ \varphi_{2}^{-1}\circ \dots \circ
\varphi_{n+1}^{-1}.$$

Then $\varphi^{-1}$
is a conformal map of $\H$ onto a region $\Om_c$ such that
$z_j\in \b\Om_c$, $j=0, \dots, n$. 
The portion $\gamma_j$ of $\b\Om_c$ between 
$z_j$ and $z_{j+1}$ is the image of the arc of a circle in the
upper half plane by the
analytic map $\f_1^{-1}\circ\dots\circ \f_j^{-1}$. 
In more picturesque language, after applying 
$\varphi_1$, we grab the ends of the displayed
horizontal line segment and pull, splitting apart  or unzipping 
the curve at $0$.
The remaining data points move down until they hit $0$ and then
each splits into two points, one on each side of $0$, moving further
apart as we continue to pull. 

As an aside, we make a few comments. As mentioned $\b\Om_c$ is piecewise
analytic. It is easy to see that it is also $C^1$ since the inverse of
the basic map $f_a$ in Figure 2 doubles angles at $0$ and halves
angles at $\pm c$. In fact it is also $C^{3\over 2}$ (see Proposition 3.12).
If the data points $\{z_j\}$ lie on the boundary of a given region
$\b\Om$, 
the analyticity of $\b\Om_c$ also allows us in many situations 
(see Proposition 2.5 and Corollary 3.9) to 
extend $\varphi_c$ analytically across $\b\Om_c$ so that the extended
map is a conformal map of $\Om$ onto a region with boundary very
close to $\b\D$. Note
also that $\varphi$ is a conformal map of the complement of $\Om_c$, 
$\C^*\setminus \overline{\Om_c}$, onto the lower half plane, $\C\setminus \overline{\H}$ 
 where $\C^*$ denotes the extended plane.  Simply follow the unshaded
region in $\H$ in Figure 3.
Finally, we remark that it is easier to use geodesic arcs in 
the right-half plane instead of in the
upper-half plane when coding the algorithm, because most computer
languages adopt the convention $-{{\pi}\over 2} < \arg \sqrt z \le 
{{\pi}\over 2}$.

\bigskip
\noindent{\bf The Slit Algorithm}
\bigskip

Given a region $\Om$, then we can select boundary points $z_0,
\dots, z_n$ on $\b\Om$ and apply the geodesic algorithm.  
We can view the circular arcs
$\gamma$ for the basic maps $f_a$ as approximating the image of the
boundary of $\Om$ between $0$ and $a$ with a circular
arc at each stage. 
We can improve the approximation 
by using  straight lines instead of orthogonal arcs.
So in the slit algorithm we replace the inverse
of the maps
$f_a$ by conformal maps  $g_a:\H\longrightarrow\H\setminus {\rm L}$
where ${\rm L}$ is
a line segment from $0$ to $a$. Explicitly

$$g_a(z) =C(z-p)^p(z+1-p)^{1-p}$$

\noindent where $p =\arg{a}/\pi$ and $C = |a|/p^p (1-p)^{1-p}$. 

\bigskip
\psfrag{H}{$\H$}
\psfrag{L}{${\rm L}$}
\psfrag{HL}{$\H\setminus {\rm L}$}
\psfrag{bm1}{$p-1$}
\psfrag{b}{$p$}
\psfrag{0}{$0$}
\psfrag{a}{$a$}
\psfrag{ga}{$g_a(z)$}
\psfrag{ppi}{$\pi p$}
\centerline{\includegraphics[width=5.0in]{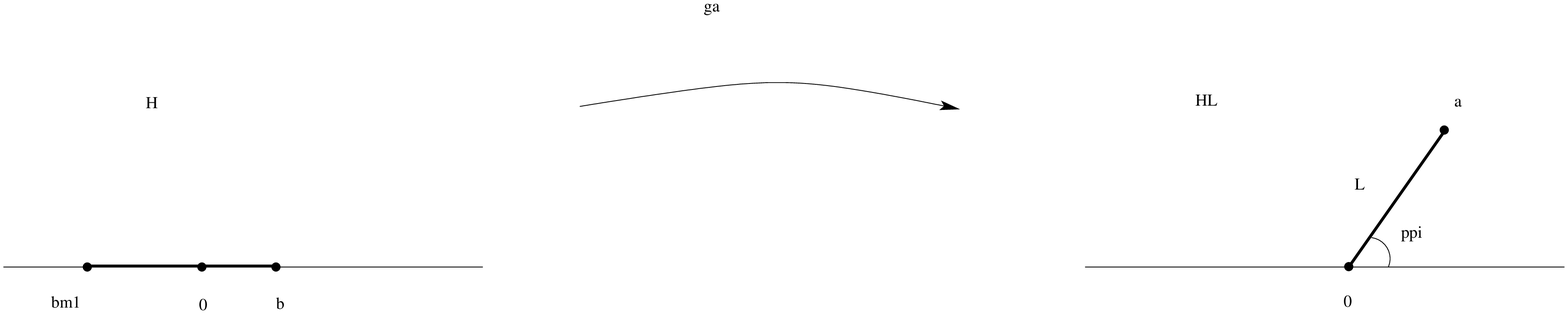}}
\nobreak
\centerline{{\bf Figure 4.} The Slit Maps.}
\vskip 0.5truein

One way to see that $g_a$ is a conformal map, is to 
note that as $x$ traces the real
line from $-\infty$ to $+\infty$, ~$g_a(x)$ traces the boundary of
$\H\setminus {\rm L}$ and $g_a(z)\sim Cz$ for large $z$ 
and then apply the argument principle. 
Another method would be to construct $Re \log g_a$ using harmonic
measure as in the first two pages of [GM].
As in the basic maps of the geodesic algorithm, 
the line segment from $0$ to $a$ is
opened to two adjacent intervals intervals on $\R$ by $f_a=g_a^{-1}$ with 
$f_a(a)=0$ and
$f_a(\infty)=\infty$.  The map $f_a$ cannot be written in terms of
elementary functions, but an effective and rapid numerical inverse will be
described in section 4.

We note that as in the geodesic algorithm, the boundary of the region
$\Om_c$ computed with the slit algorithm will be piecewise analytic.
However it will not be $C^1$. A curve is called $C^1$ if the arc
length parameterization has a continuous first derivative. In other
words, the
direction of the unit tangent vector is continuous. 
Indeed if $g_a$ is the map illustrated
by Figure 4, and if $g_b$ is another such map then $g_b\circ g_a$
forms a curve with angles $2\pi p$ and $2\pi (1-p)$ on either side of
the curve at $b=g_b(0)$. Since analytic maps preserve angles, the
boundary of the computed region consists of analytic arcs with
endpoints at the data points, and angles determined by the basic maps.
This will allow us to accurately compute conformal maps to regions with 
(a finite number of) ``corners'', or ``bends''.

\bigskip 
\bigskip
\noindent {\bf The Zipper Algorithm}

\bigskip
We can further improve the approximation by replacing the linear slits
with arcs of (non-orthogonal) circles. In this version we assume there
are an even number of boundary points, $z_0, z_1, \dots, z_{2n+1}$.
The first map is replaced by
$$\varphi_1(z)=
\sqrt{{{(z-z_2)(z_1-z_0)}\over{(z-z_0)(z_1-z_2)}}}$$
which maps the complement in the extended plane 
of the circular arc through $z_0, z_1, z_2$
onto $\H$.
At each subsequent stage, instead of
pulling down one point $\zeta_k$, we can find a unique circular arc
through $0$ and the (images of) the next two data points $\zeta_{2k-1}$
and $\zeta_{2k}$. By a linear
fractional transformation ${\ell}_a$ which preserves $\H$, 
this arc is mapped to a line segment (assuming the arc is not tangent
to $\R$ at $0$. See Figure~5.

\bigskip
\psfrag{a}{$d$}
\psfrag{b}{$b$}
\psfrag{c}{$c$}
\psfrag{d}{$a$}
\psfrag{ha}{${\ell}_a(z)=\displaystyle{{z}\over{1-z/b}}$}
\psfrag{0}{$0$}
\psfrag{pip}{$\pi p$}
\psfrag{p}{$p$}
\psfrag{pm}{$p-1$}
\psfrag{gdi}{$g_d^{-1}(z)$}
\centerline{\includegraphics[width=4.8in]{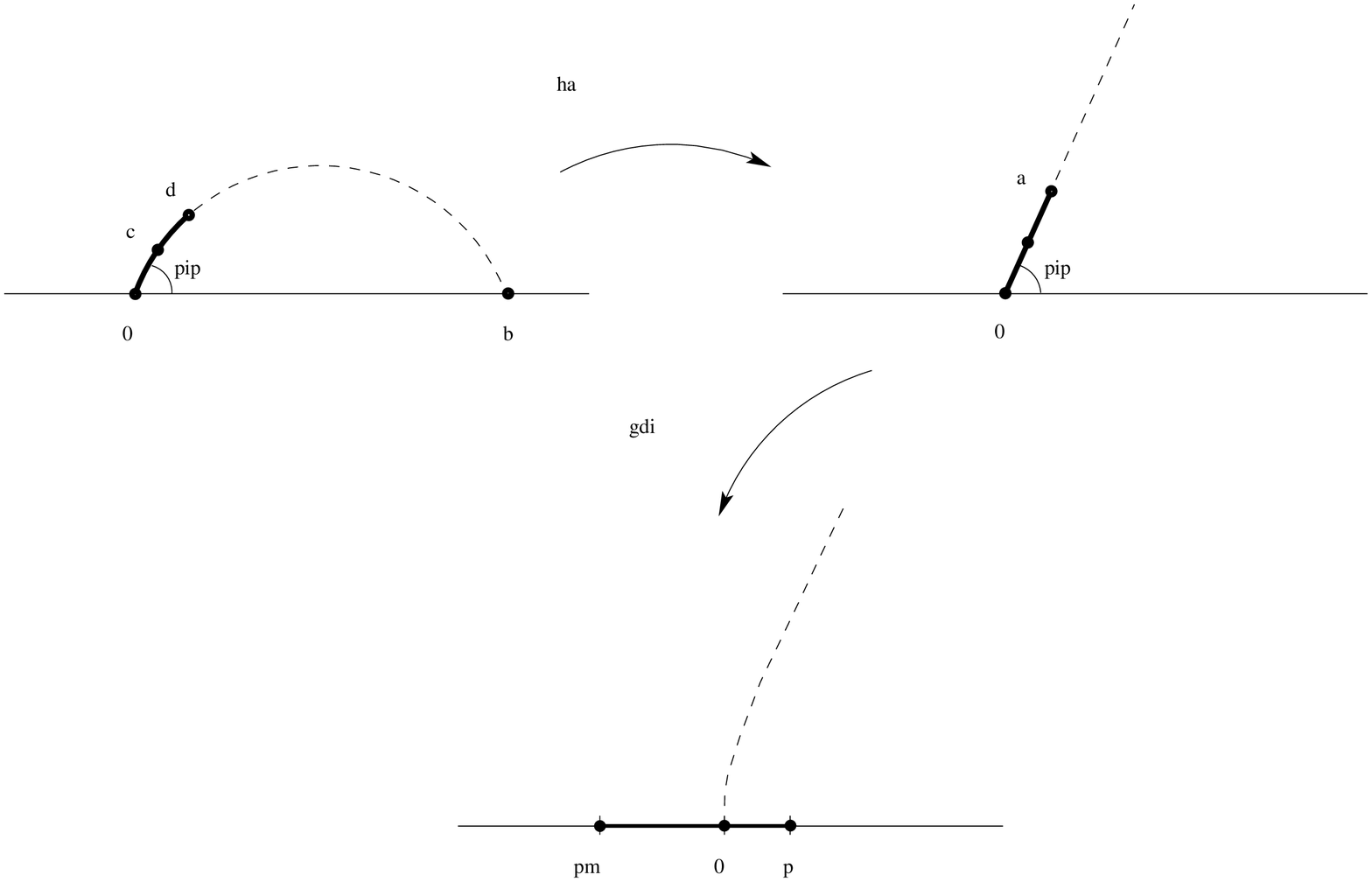}}
\nobreak
\centerline{{\bf Figure 5.} The Circular Slit Maps.}
\vskip 0.5truein

The complement of this segment in $\H$ can then be mapped to $\H$ as
described in the slit algorithm, using $g_d^{-1}$ where $d=a/(1-a/b)$.
The composition $h_{a,c}=g_d^{-1}\circ{\ell}_a$ then maps the complement
of the circular arc in $\H$ onto $\H$. Thus at each stage we are giving a
``quadratic approximation'' instead of a linear approximation to the 
(image of) the boundary. The last map $\varphi_{n+1}$ is a conformal
map of  the intersection of a disc with $\H$ where the boundary
circular arc passes through $0$, the image of $z_{2n+1}$ and the image
of $z_0$ by the composition $\varphi_n\circ\dots\circ\varphi_1$.
See Figure 6.

\vskip 0.5truein
\psfrag{z0}{$z_0$}
\psfrag{z1}{$z_1$}
\psfrag{z2}{$z_2$}
\psfrag{z3}{$z_3$}
\psfrag{z4}{$z_4$}
\psfrag{zn}{$z_n$}
\psfrag{ze3}{$\zeta_3$}
\psfrag{ze4}{$\zeta_4$}
\psfrag{ze5}{$\zeta_5$}
\psfrag{ze6}{$\zeta_6$}
\psfrag{ze7}{$\zeta_7$}
\psfrag{ze8}{$\zeta_8$}
\psfrag{p1}{$\varphi_1=\displaystyle{\sqrt{{(z-z_2)(z_1-z_0)}\over{(z-z_0)(z_1-z_2)}}}$}
\psfrag{p2}{$\varphi_2=h_{\zeta_3,\zeta_4}$}
\psfrag{p3}{$\varphi_3=h_{\zeta_5,\zeta_6}$}
\psfrag{pn}{$\varphi_n=h_{\zeta_{2n-1},\zeta_{2n}}$}
\psfrag{pnp}{$\varphi_{n+1}=C\displaystyle{\biggl({z\over{1-z/\zeta_{2n+2}}}\biggr)^{{\pi}\over{\alpha}}}$}
\psfrag{H}{$\H$}
\psfrag{al}{$\alpha$}
\psfrag{0}{$0$}
\psfrag{Om}{$\Omega_c$}
\psfrag{zep}{$\zeta_{2n+1}$}
\psfrag{ze2p}{$\zeta_{2n+2}$}
\centerline{\includegraphics[height=4.5in]{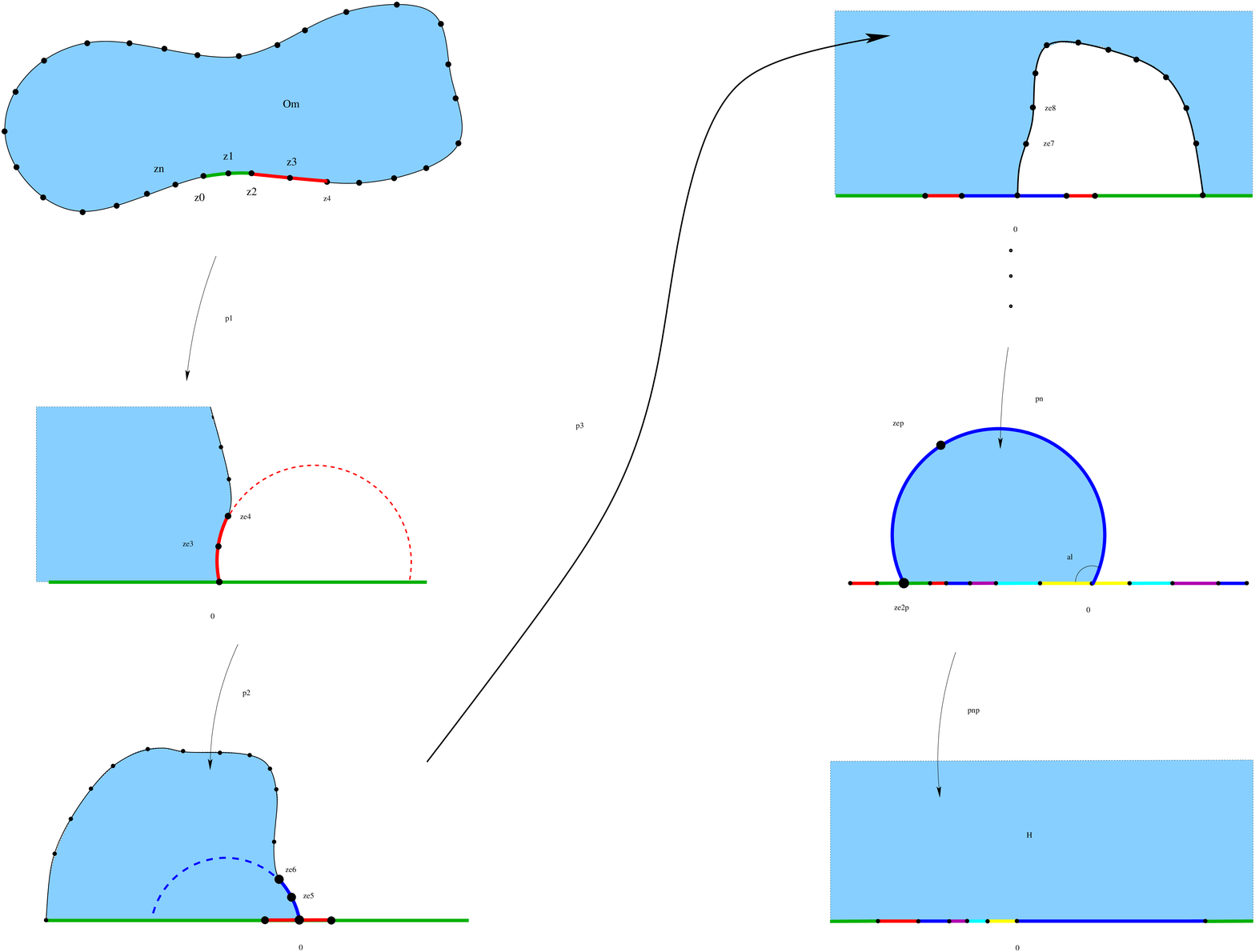}}
\nobreak
\centerline{{\bf Figure 6.} The Zipper Algorithm.}
\bigskip

If the zipper algorithm is used to approximate the boundary of a
region with bends or angles at some boundary points, then better
accuracy is obtained if the bends occur at even numbered vertices
$\{z_{2n}\}$.
\bigskip

\bigskip 
\bigskip
\noindent {\bf Conformal Welding}
\bigskip 

The discovery of the slit algorithm by the first author came from 
considering conformal weldings.  (The simpler geodesic
algorithm was discovered later.) A
decreasing continuous function $h:[0,+\infty) \to (-\infty,0]$ with
$h(0)=0$ is called a {\bf conformal welding} if there is a conformal
map $f$ of $\H$ onto $\C\setminus \gamma$ where $\gamma$ is a Jordan
arc from $0$ to $\infty$ such that $f(x)=f(h(x))$ for $x\in \R$. 
In other words, the
map $f$ pastes the negative and positive real half-lines together
according to the prescription $h$ to form a curve. One way to
approximate a conformal welding is to prescribe the map $h$ at
finitely many points and then construct a conformal mapping of $\H$
which identifies the associated intervals. 

A related problem, which
the first author considered in joint work with L. Carleson, is: given 
angles $\alpha_1, \alpha_2,\dots,\alpha_n$ and 
$0 < x_1 < x_2 < \dots < x_n$, find points $y_n < \dots <  y_1 < 0$ so
that there is a Schwarz-Christoffel map $f$ of $\H$ onto a region bounded
by a polygonal arc tending to $\infty$ with angles $\alpha_j,
2\pi-\alpha_j$ at the $jth$ vertex $f(x_j)=f(y_j)$. This map welds the
intervals $[x_j,x_{j+1}]$ and $[y_{j+1},y_j]$, $j=1,\dots,n$. 
Unfortunately, at the time the best Schwarz-Christoffel method was only fast
enough to do this problem with polygonal curves with up to 20 bends.

The basic maps $g_a$ can be used to compute the conformal maps of weldings.
Indeed, suppose $y_1<0<x_1$, let $a=x_1/(x_1-y_1)$, and apply the map
$g_a(z/(x_1-y_1))$. This map identifies the intervals $[y_1,0]$ and
$[0,x_1]$, by mapping them to the two ``sides'' of a line segment
$L\subset\H$. Composing maps of this form will give a conformal map
$\varphi:\H \to \C\setminus \gamma$ such that 
$\varphi([x_j,x_{j+1}])=\varphi([y_{j+1},y_j]).$ The final intervals
are welded together using the map $z^2$.
The numerical computation of these maps is easily fast enough to
compose $10^5$  basic maps, thereby giving an approximation to almost
any conformal welding. Conversely, given a Jordan arc $\gamma$
connecting $0$ to $\infty$, the
associated welding can be found approximately by
using the slit algorithm to approximate the
conformal map from $\H$ to the complement of $\gamma$.

The idea of closing up such a region using a map of the form
$\varphi_{n+1}$ was suggested by L. Carleson, for which we thank him.

Since we have been asked about this a couple of times, we note that
conformal welding can also be defined using the conformal maps 
to the
inside and outside of a closed Jordan curve. If $f$ is a conformal map of
the unit disc
$\D$ onto a Jordan region $\Om$ and if $g$ is a conformal map of
$\C\setminus\overline{\D}$ onto $\C\setminus\overline{\Om}$ which
maps $\infty$ to $\infty$ then $f$ and $g$ extend to be homeomorphisms
of $\overline{\D}$ onto $\overline{\Om}$ and $\C\setminus \D$ onto 
$\C\setminus \Om$
respectively. Then the map 
$$h= f^{-1} \circ g :\b \D \longrightarrow \b \D$$
is a homeomorphism of the unit circle and is also called a conformal
welding. Again, if we approximate a homeomorphism $h$ by prescribing
it at finitely many points on the circle, then we can use the slit
algorithm to identify the corresponding intervals. Simply
map the disc to the upper-half plane so that the first interval $I$ is
mapped to $\R^+$, the positive reals, and map the complement of the disc to
the lower half plane so that desired image $h(I)$ is mapped to
$\R^+$. Apply $i\sqrt{z}$ and now proceed to identify the remaining
intervals as above. Conformal welding can also be accomplished using
the geodesic algorithm. We leave the elementary details to the
interested reader.

\bigskip

From this point of view, the slit or the geodesic algorithms find the
conformal welding of a curve (approximately). From the point of view
of increasing the boundary via a small curve $\gamma_j$ from $z_j$ to
$z_{j+1}$, the algorithms
are discrete solutions of L\"owner's differential equation.

\bigskip
\bigskip
\S {\bf 2. Disc-chains}
\bigskip

The geodesic algorithm can be applied to any sequence of data points
$z_0, z_1, \dots, z_n$, unless the points are out of order in the
sense that a data point $z_j$ belongs to the
geodesic from $z_{k-1}$ to $z_k$, for some $k < j$. 
In this section we will give a simple condition on the data points
$z_0, z_1,\dots,z_n$ which is sufficient to guarantee that the
curve computed by the geodesic algorithm 
is close to the polygon with vertices $\{z_j\}$.

\bigskip
\proclaim Definition 2.1. A {\bf disc-chain} $D_0, D_1,\dots,D_n$ is a
sequence of pairwise disjoint open discs such that 
$\b D_j$ is tangent to $\b D_{j+1}$, for $j=0,\dots,n-1$. 
A {\bf closed disc-chain} is a disc-chain
such that $\b D_n$ is tangent to $\b D_0$.
\bigskip

Any closed Jordan polygon $P$, for example, can be covered by a closed
disc-chain
with arbitrarily small radii and centers on $P$. There are several ways to
accomplish this, but one straightforward method is the following:
Given $\eps >0$, 
find pairwise disjoint discs $\{B_j\}$ centered at each vertex, and
of radius less than $\eps$. Then 
$$P\setminus \bigcup_j B_j = \bigcup L_k$$
where $\{L_k\}$ are pairwise disjoint closed line segments. Cover each
$L_k$ with a disc-chain centered on $L_k$ tangent to the
corresponding $B_j$ at the ends, and radius less than half the
distance to any other $L_i$, and less than $\eps$.

\vskip 0.5truein
\psfrag{}{}
\centerline{\includegraphics[height=1.5in]{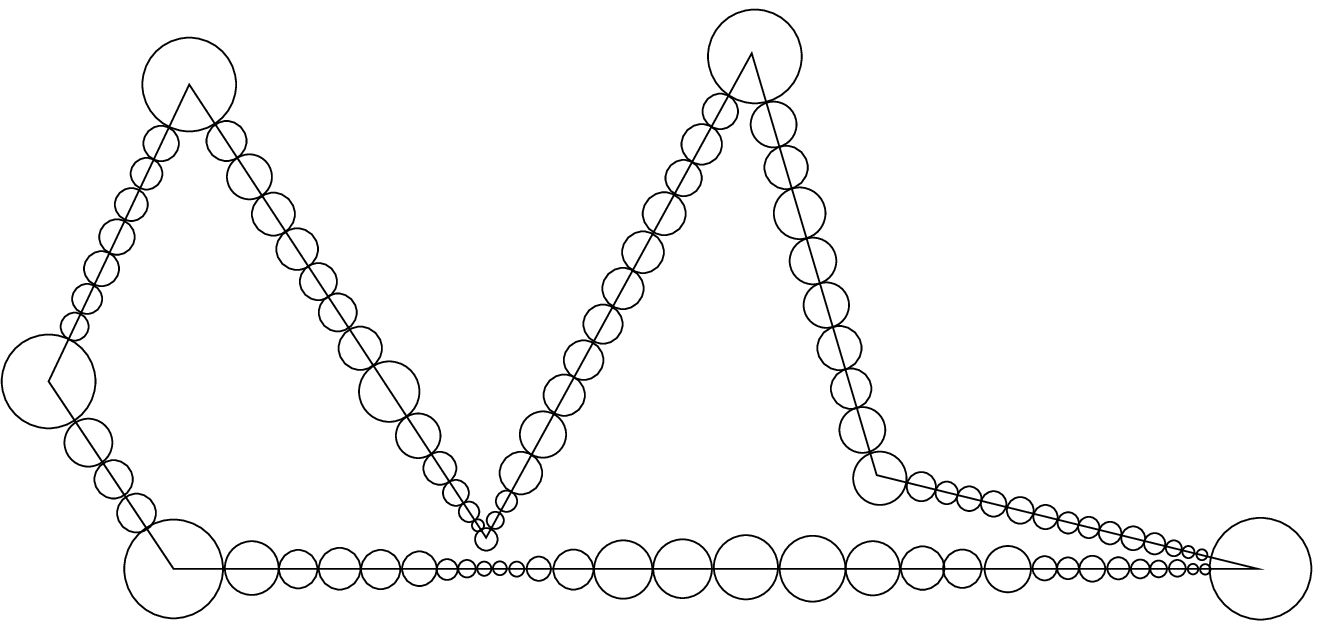}}
\nobreak
\centerline{{\bf Figure 7.} Disc-chain covering a polygon.}
\vskip 0.5truein

Another method for constructing a disc-chain is to use a Whitney
decomposition of a simply connected domain. Suppose $\Om$ is a 
simply connected domain contained in the unit square. The square is
subdivided into 4 equal squares. Each of these squares is subdivided
again into 4 equal squares, and the process is repeated. If $Q$ is a
square, let $2Q$ denote the square with the same center, and sides
twice as long. 
In the subdivision process, if a square $Q$ satisfies $2Q\subset\Om$, 
then no further subdivisions are made in $Q$. Let $U_n$ be the
union of all squares $Q$ obtained by this process
with side length at least $2^{-n}$ for which $2Q\subset \Om$.
If $z_0\in \Om$, let $\Om_n$ be the component of the interior of $U_n$
containing $z_0$. Then $\b\Om_n$ is a polygonal Jordan
curve. Note that $\b\Om_n$ consists of sides of squares $Q$ with
length $2^{-n}$. Thus we can form a disc chain by placing a disc of
radius $2^{-n}/2$ at each vertex of $\b\Om_n$. The points of tangency
are the midpoints of each square with edge length $2^{-n}$ on
$\b\Om_n$.

Yet another method for constructing a disc-chain would be to start with a
hexagonal grid of tangent discs, all of the same size, then select a
sequence of these discs which form a disc-chain. The boundary circles
of a circle packing of a simply connected domain can also be used to
form a disc-chain. See for example any of the pictures in Stephenson
[SK].



\bigskip

If $D_0, D_1,\dots,D_n$ is a closed disc-chain, set
$$z_j=\b D_j \cap \b D_{j+1},$$ 
for $j=0,\dots,n$, where $D_{n+1}\equiv D_0$.

\bigskip
\proclaim Theorem 2.2.  If  $D_0, D_1,\dots,D_n$ is a closed disc-chain, 
then the geodesic algorithm applied to the data $z_0, z_1, \dots,z_n$
produces a conformal map $\varphi_c^{-1}$ from the upper half plane $\H$
to a region bounded by a $C^1$ and piecewise analytic
Jordan curve $\gamma$ with 
$$\gamma \subset \bigcup_0^n (D_j \cup z_j).$$ 
\bigskip

\proof An arc of a circle which is orthogonal to $\R$ is a hyperbolic
geodesic in the upper half plane $\H$. Let $\gamma_j$ denote the
portion of the computed boundary, $\b\Om_c$, between $z_j$ and $z_{j+1}$.
Since hyperbolic geodesics are preserved by conformal maps,
$\gamma_j$ is a hyperbolic geodesic in 
$$\C^*\setminus \cup_{k=0}^{j-1} \gamma_k.$$
For this reason, we call the algorithm the ``geodesic'' algorithm.

Using the notation of Figure 2, 
each map $f_a^{-1}$ is analytic across $\R\setminus
\{\pm c\}$, where $f_a(\pm c)=0$,
and $f_a^{-1}$ is approximated by a square
root near $\pm c$. If $f_b^{-1}$ is another
basic map, then $f_b^{-1}$ is analytic and asymptotic to a
multiple of $z^2$ near $0$. Thus $f_b^{-1}\circ
f_a^{-1}$ preserves angles at $\pm c$. The geodesic $\gamma_j$ then is
an analytic arc which 
meets $\gamma_{j-1}$ at $z_j$  with angle $\pi$. Thus 
the computed boundary $\b\Om$ is $C^1$ and piecewise analytic. 
The first arc
$\gamma_0$ is a chord of $D_0$ and hence not tangent to $\b D_0$. 
Since the angle at $z_1$ between $\gamma_0$ and
$\gamma_1$ is $\pi$, ~$\gamma_1$ must enter $D_1$,
and so by J\o rgensen's theorem (see Theorem
A.1 in the appendix) 
$$\gamma_1 \subset D_1,$$
and $\gamma_1$ is not tangent to $\b D_1$. By induction
$$\gamma_j\subset D_j,$$
$j=0,1,\dots,n$.  
\endproof

\bigskip Disc-chains can be used to approximate the boundary of an
arbitrary simply connected domain.

\bigskip
\proclaim Lemma 2.3. Suppose that  $\Om$ is a bounded simply connected domain.
If $\eps >0$, then there is a
disc-chain $D_0,\dots,D_n$ so that the radius of each $D_j$ 
is at most $\eps$ and
$\b \Om$ is contained in an $\eps$-neighborhood of $\cup D_j$.
\bigskip

\proof 
We may suppose that $\Om$ is contained in the unit square. Then for $n$
sufficiently large, the disc chain constructed using the Whitney
squares with side length at least $2^{-n}$, as described above,
satisfies the conclusions of Lemma 2.3.
\endproof 
\bigskip

The {\bf Hausdorff distance} $d_H$ in a metric $\rho$ between two sets $A$ and $B$ is the smallest number $d$
such that every point of $A$ is within $\rho$-distance $d$ of $B$, and every point of $B$ is within
$\rho$-distance $d$ of $A$. The $\rho$-metrics we will consider in this article are the Euclidean
and spherical metrics.

A consequence 
is the following theorem.

\bigskip
\proclaim Theorem 2.4. If $\Om$ is a bounded simply connected domain 
then for any
$\eps>0$, the geodesic algorithm can be used to find a conformal map
$f_c$ of $\D$ 
onto a Jordan region $\Om_c$ so that 
$$d_H(\b\Om,\b\Om_c) < \eps,\eqno(2.1)$$
where $d_H$ is the Hausdorff distance in the Euclidean metric.
If $\b\Om$ is a Jordan curve then we can find $f_c$ so that
$$\sup_{z\in\overline \D} |f(z)-f_c(z)| < \eps,$$
where $f$ is a conformal map of $\D$ onto $\Om$.
\bigskip

\proof The first statement follows immediately from Theorem 2.2 and
Lemma 2.3. To prove the second statement, note that the boundary of
the regions constructed with the Whitney decomposition converges to
$\b\Om$ in the Fr\'echet sense. By a theorem of Courant [T, page 383], 
the mapping
functions can be chosen to be uniformly close. \endproof

\bigskip
We note that if $\Om$ is unbounded, Lemma 2.3 and Theorem 2.4 remain 
true if we use
the spherical metric instead of the Euclidean metric to measure the
radii of the discs and the distance to $\b\Om$.
\bigskip
There are other ways besides using the Whitney
decomposition to approximate the boundary of a region by a disc-chain
and hence to approximate the mapping function. However, Theorem 2.4
does not give an explicit estimate of the distance between mapping
functions in terms of the geometry of the regions. This issue will be
explored in greater detail in Sections 5 and 6.

\bigskip
The von Koch snowflake is an example of a simply connected Jordan domain
whose boundary has Hausdorff dimension $> 1$.
The standard construction of the von Koch snowflake
provides a sequence of polygons which approximate it. By Theorem 2.4
the mapping functions constructed from these disc-chains converge
uniformly to the conformal map to the snowflake.

\vskip 0.5truein
\psfrag{}{}
\centerline{\includegraphics[height=1.5in]{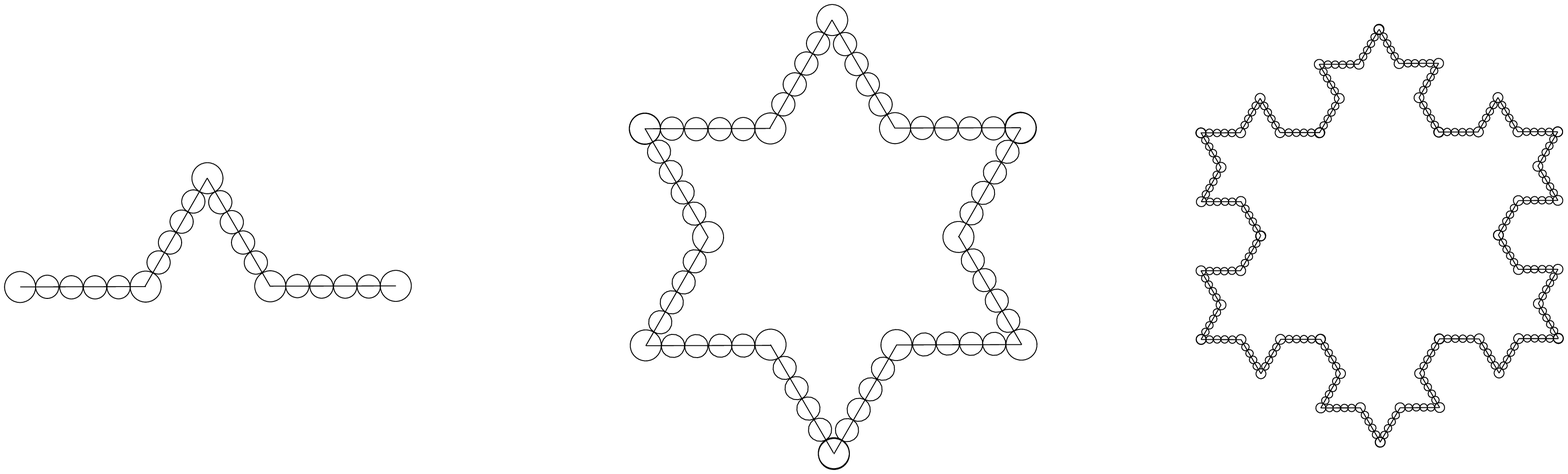}}
\nobreak
\centerline{{\bf Figure 8.} Approximating the von Koch snowflake.}
\vskip 0.5truein

It is somewhat amusing and perhaps
known that a
 constructive proof of the Riemann mapping theorem (without the use
of normal families)
then follows. 
Using linear fractional transformations and a square
root map, we may suppose $\Om$ is a bounded simply connected domain. 
Using the disc-chains associated with increasing levels of the 
Whitney decomposition for instance, 
$\Om$ can be exhausted by an increasing sequence of
domains $\Om_n$ for which the geodesic algorithm can be used to compute the
conformal map $f_n$ of $\Om_n$ onto $\D$ with
$f_n(z_0)=0$ and $f_n'(z_0)>0$.  Then by Schwarz's lemma
$$u_n(w)=\log\biggl|{f_m(w)\over{f_n(w)}}\biggl|$$
for $n=m+1,m+2,\dots$ is an increasing sequence of positive harmonic 
functions on $\Om_m$ which is bounded above at $z_0$ by 
Schwarz's lemma applied to $f_n^{-1}$, since $\Om$ is bounded. 
By Harnack's estimate $u_n$ is bounded on compact
subsets of $\Om$ and by the Herglotz integral formula, 
$\log{f_m(w)\over{f_n(w)}}$ converges uniformly
on closed discs contained in $\Om_m$. Thus $f_n$ converges
uniformly on compact subsets of $\Om$ to an analytic function $f$.
By Hurwitz's theorem  $f$ is one-to-one and by Schwarz's lemma
applied to $f_n^{-1}$, $f$ maps $\Om$ onto $\D$.

In the geodesic algorithm, we have viewed the maps $\varphi_c$ and 
$\varphi_c^{-1}$ as conformal maps between $\H$ and a region $\Om_c$ whose
boundary contains the data points. If we are given a region $\Om$, and
choose data points $\{z_k\}\in \b\Om$ properly, then the next
proposition says that the computed
maps $\varphi_c$ and $\varphi_c^{-1}$ are also  conformal maps between 
the original region $\Om$ and a region ``close'' to $\H$.

\bigskip
\proclaim Proposition 2.5.  If $D_0,\dots,D_n$ 
is a closed disc-chain with points of tangency $\{z_k\}$,
and if $\Om$ is a simply connected domain such that
$$\b\Om\subset \bigcup_{k=0}^{n} \overline{D_k}$$
then the computed map $\varphi_c$ for the data points $\{z_k\}_0^{n}$
extends to be conformal on $\Om$. 
\bigskip

We remark that changing the sign of the last map $\varphi_{n+1}$ in
the construction of $\varphi_c$ gives a conformal map of the complement
of the computed region onto $\H$. We choose the sign so that the
computed boundary winds once around a given interior point of $\Om$.

\proof
Without loss of generality $\Om\supset \bigcup_{k=0}^{n}
D_k$ and hence $\b\Om \subset\cup \b D_k$. The basic map $f_a$ in Figure 2 
extends by reflection to be a conformal
map of $\C^*\setminus (\gamma\cup \gamma^R)$ onto $\C^*\setminus
[-c,c]$, where $\gamma^R$ is the reflection of $\gamma$ about $\R$.
In fact, if $\sigma\subset \H$ 
is any connected set such that $0$, $a \in \overline \sigma$ then $f_a$ is
conformal on $\C^*\setminus (\sigma\cup \sigma^R)$, where $\sigma^R$ is the
reflection of $\sigma$ about $\R$. In particular, if $U$ is a simply
connected region contained in $\H$ with $0,a\in \b U$, then
$\C\setminus \overline{f_a(U\cup U^R)} $ consists of two open sets $V \cup -V$
where $(0,c)\in V$ and $(0,-c) \in -V$. Here $U^R$ denotes the
reflection of the set $U$ about $\R$ and $-V=\{-z:z\in V\}$.

Set $$\psi_k\equiv\varphi_k\circ\cdots \circ \varphi_1$$ 
and 
$$W_k=\psi_k(\C^*\setminus\{D_0\cup \dots\cup D_n\}).$$
Then we claim $\C^*\setminus\{W_k\cup W_k^R\}$ consists of $2(n+1)$
pairwise disjoint simply connected regions:
$$\C^*\setminus\{W_k\cup W_k^R\}=\bigcup_{j=k}^n \psi_k(D_j) \cup \psi_k(D^j)^R
\cup \bigcup_{j=1}^{2k} U_{k,j},$$
where each region $U_{k,j}$ is symmetric about $\R$ and $\R\subset
\cup_{j=1}^{2k} \overline{U_{k,j}}$. The case $k=1$ follows since
$\psi_1(\C^*\setminus D_0)$ is bounded by two lines from $0$ to
$\infty$. As noted above, the image of $\psi_{k-1}(D_{k-1})\cup
\psi_{k-1}(D_{k-1})^R$  by the map
$\varphi_k$ consists of two regions $V$ and $-V$. The claim now
follows by induction.
Note $\delta_k=\psi_k(\b\Om\cap \b D_{k-1})$ 
is a simple curve connecting a point
$-c_k$ to $0$ and $0$ to $c_k$, satisfying $z\in \delta_k$ if and only
if $-z\in \delta_k$, since $\sqrt{z^2+1}$ is odd.  Since  each
$\varphi_k$ extends to be one-to-one and analytic on
$\psi_{k-1}(D_{k-1})$
and since $U_{n,j}$, $j=1,...2n$ are disjoint, the map $\psi_n$ is
one-to-one and analytic on $\Om$. By direct inspection, the final map
$\varphi_{n+1}$ extends to be one-to-one and analytic, completing the
proof of Proposition 2.5 
\endproof
\bigskip
As one might surmise from the proof of Proposition 2.5, care must be
taken in any numerical implementation to assure that the proper
branch of $\sqrt{z^2+c^2}$ is chosen at each stage  in order to find
the analytic extension of the computed map to all
of $\Om$.

\bigskip
\bigskip
\S {\bf 3. Diamond-chains and Pacmen}
\bigskip

If we have more control than the disc-chain condition on the behavior
of the boundary of a region, then we show in this section that the geodesic
algorithm approximates the boundary with better estimates. 
We will first restrict our attention to domains of the form 
$\C\setminus \gamma$ where $\gamma$ is a Jordan arc tending to
$\infty$.

\bigskip
\proclaim Definition 3.1. An {\bf $\eps$-diamond} $D(a,b)$ is an open
rhombus with opposite
vertices $a$ and $b$ and interior angle $2\eps$ at $a$ and at $b$. 
If $a=\infty$, then an {\bf $\eps$-diamond} $D(\infty,b)$  is a sector
$\{z:|\arg (z-b) -\theta|<\eps\}$. An
{\bf $\eps$-diamond-chain} is a pairwise disjoint sequence of
$\eps$-diamonds $D(z_0,z_1), D(z_1,z_2), \dots D(z_{n-1},z_n)$.
A {\bf closed $\eps$-diamond-chain} is an $\eps$-diamond-chain with
$z_{n}=z_0$. 
\bigskip

See Figure 9. Let $B(z,R)$ denote the disc centered at $z$ with radius $R$.

\bigskip
\proclaim Definition 3.2. A {\bf pacman} is a region of the form
$$P=B(z_0,R)\setminus\{z:|\arg(\overline{\lambda}(z-z_0))|\le\eps\},$$
for some {\bf radius} $R < \infty$, {\bf center} $z_0$, 
{\bf opening} $2\eps>0$, 
and {\bf rotation} $\lambda$, $|\lambda|=1$.
\bigskip

Let $C_1$ be a constant to be chosen later (see Lemma 3.7),
and let $z_0=\infty.$

\bigskip
\proclaim Definition 3.3. We say that an $\eps$-diamond-chain
$D(\infty,z_1), D(z_1,z_2), \dots D(z_{n-1},z_n)$,
satisfies the {\bf $\eps$-pacman condition} if for
each $1\leq k\leq n-1$ the pacman 
$$P_k =
B(z_k,R_k)\setminus\{z:|\arg\Bigl({{z-z_k}\over{z_k-z_{k+1}}}\Bigr)|\le
{{\eps}}\},$$ 
with radius $R_k=C_1|z_{k+1}-z_k|/\eps^2$ satisfies
$$\biggl(\bigcup_{j=0}^{k-2} D(z_j,z_{j+1})\biggr) \cap P_k
=\emptyset.$$
\bigskip

The pacman $P_k$ in Definition 3.3 is chosen to be symmetric about the
segment between $z_k$ and $z_{k+1}$ 
with opening $2\eps$ equal to the interior angle $2\eps$ in the
diamond-chain.
Note that the $\eps$-diamond
$D(z_{k-1},z_k)$ may intersect $P_k$. 
%

\vskip 0.5truein
\psfrag{zk}{$z_k$}
\psfrag{zkm}{$z_{k-1}$}
\psfrag{zkp}{$z_{k+1}$}
\psfrag{D01}{$D(z_0,z_1)$}
\psfrag{z1}{$z_0$}
\psfrag{z0}{$z_1$}
\psfrag{Pk}{$P_k$}
\psfrag{Rk}{$R_k$}
\centerline{\includegraphics[height=2.0in]{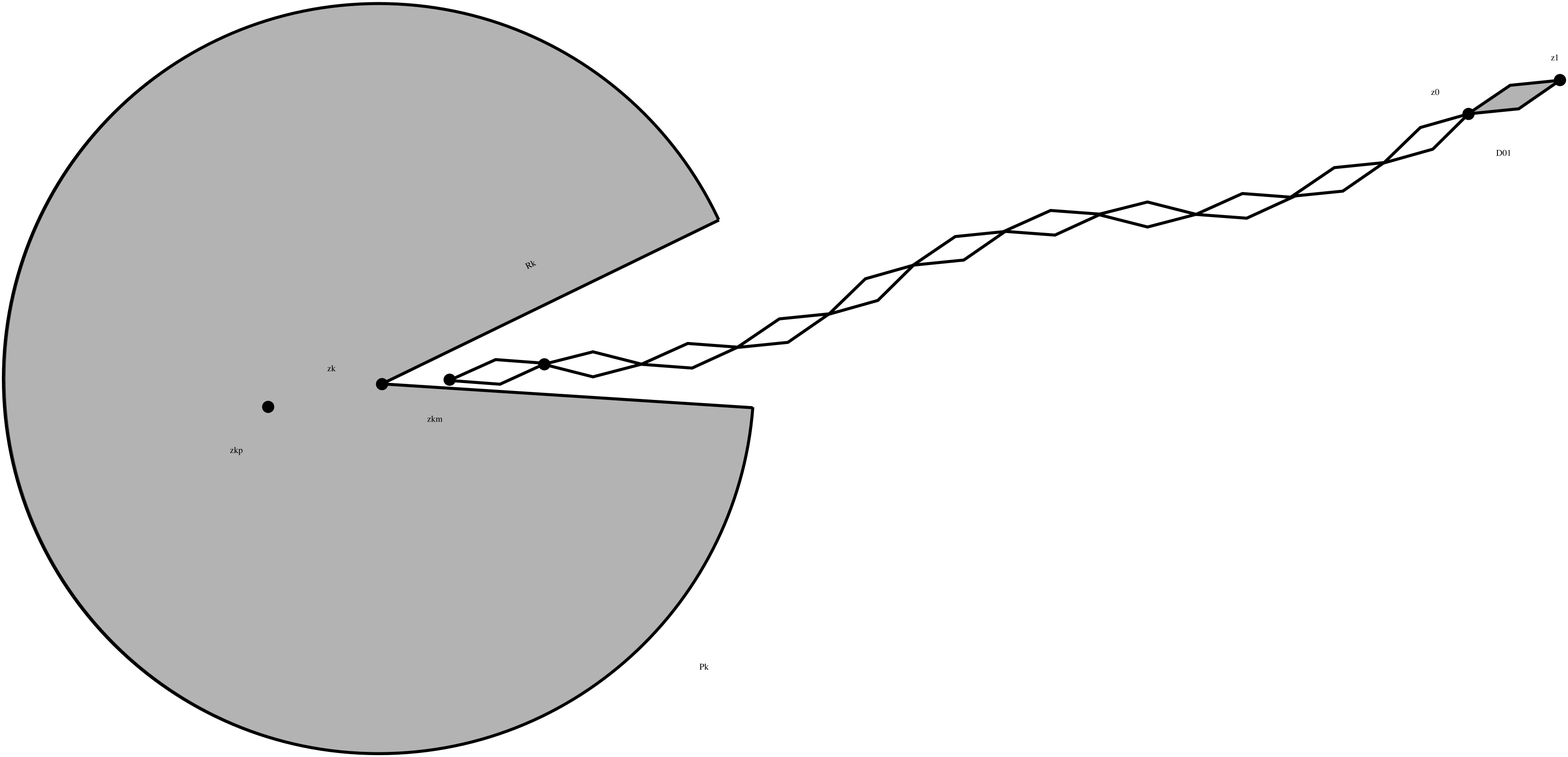}}
\nobreak
\centerline{{\bf Figure 9.} A Diamond-chain and a Pacman.}
\vskip 0.5truein

When $z_0=\infty$, the first map in the geodesic algorithm
is replaced by $\varphi_1(z)=\lambda\sqrt{z-z_1}$. 
The argument of $\lambda$ can be
chosen so that $\varphi_1(z_2)$ is purely imaginary, in which case the
boundary of the constructed region contains the half-line from $z_2$
through $z_1$ and $\infty$.  We will henceforth assume that
$$D(\infty, z_1)=\{z: |\arg\biggl({{z-z_1}\over{z_1-z_2}}\biggr)|<
\eps\}.$$

\bigskip
\proclaim Theorem 3.4.  There exist universal constants $\eps_0>0$ 
and $C_1$ such that if an $\eps$-diamond-chain
$$D(\infty,z_1), D(z_1,z_2), \dots ,D(z_{n-1},z_n)$$
satisfies the $\eps$-pacman condition with 
$\eps < \eps_0$, and if 
$$\Biggl|\arg\biggl({{z_{k+1}-z_k}\over{z_k-z_{k-1}}}\biggr)\Biggl| <
{{\eps}\over{10}},\eqno(3.1)$$
for $k=2,\dots,n-1$,
then the boundary curve $\gamma$
computed by the geodesic algorithm with the data $z_0=\infty,z_1,\dots,z_n$
satisfies
$$\gamma\subset\bigcup_{k=1}^n
\biggl(D(z_{k-1},z_k)\cup\{z_k\}\biggr).$$
Moreover, the argument $\theta$ of the tangent to $\gamma$ between 
$z_k$ and $z_{k+1}$ satisfies $|\theta - \arg(z_{k+1}-z_k)| < 3\epsilon.$
\bigskip

To prove Theorem 3.4, we first give several lemmas.

\bigskip
\proclaim Lemma 3.5. There exists $\eps_0 > 0$ such that if $\eps <
\eps_0$, and if $\Om$ is a simply connected region
bounded by a Jordan arc $\b\Om$ from $0$ to $\infty$ with
$$\{z: |\arg z| < \pi - \eps\} \subset \Om,$$
then the conformal map $f$ of $\H^+=\{z:\Re z > 0\}$ onto $\Om$
normalized so that $f(0)=0$ and $f(\infty)=\infty$ satisfies
$$|\arg z_0^2 f'(z_0)| < {{5\eps}\over 6},\eqno(3.2)$$
where $z_0= f^{-1}(1)$.
\bigskip

The circle $C_{z_0}$ which is orthogonal to the imaginary axis at
$0$ and passes through $z_0$ has a tangent vector at $z_0$ with
argument equal to $2\arg z_0$.
The quantity $\arg z_0^2 f'(z_0)$ in (3.2) 
is the argument of the tangent vector to $f(C_{z_0})$ at $f(z_0)$. 
\bigskip

\proof We may suppose that $|z_0|=1$.  Set 
$$g(z)= \log{{f(z)}\over {z^2}}.$$
Then $|\Im g(z)| \le \eps$ on $\b\H^{+}$ and hence also on $\H^{+}$,
and $|\arg{z_0}| \le {{\eps}\over {2}},$ since $f(z_0)=1$. Set
$\alpha={{\pi}\over {2\eps}}$ and
$$A =e^{\alpha g(z_0)}=z_0^{-2\alpha},$$
$$\varphi(z)= {{e^{\alpha z} - A}\over
{e^{\alpha z} + \overline{A}}},$$
and
$$\tau(z)= {{1+z}\over{1-z}}\Re z_0 + i \Im z_0.$$

\noindent Then $\tau$ is a conformal map of $\D$ onto $\H^+$ such that
$\tau(0)=z_0$ and $\varphi$ is a conformal map of the strip
$\{|\Im z| < \eps\}$ onto $\D$ so that $\varphi(g(z_0))=0$. Thus
$h=\varphi\circ g \circ \tau$  is analytic on $\D$, bounded by $1$ and
$h(0)=0$, so that by Schwarz's lemma

$$|\varphi'(g(z_0))||g'(z_0)||\tau'(0)|=|h'(0)|\le 1.$$
Consequently
$$\biggl|{{f'(z_0)}\over{f(z_0)}} - {2\over{z_0}}\biggr| = |g'(z_0)|
\le {{2\eps |\Re A|}\over{\pi \Re z_0}}\le {{2\eps}\over{\pi
\cos{{{\eps}\over 2}}}},$$

and hence
$$\eqalign{|\arg z_0^2 f'(z_0)| &= \biggl|\arg {z_0} + \arg{{z_0
f'(z_0)}\over{f(z_0)}}\biggr|\cr
&\le {{\eps}\over 2} +
\sin^{-1}\biggl({{\eps}\over{\pi \cos{{\eps}\over 2}}}\biggr)\cr 
&= \biggl({1\over 2}+ {1\over \pi}\biggr)\eps + {\rm O}(\eps^2).}$$
This proves Lemma 3.5 if $\eps$ is sufficiently small. 
\endproof

\bigskip
\proclaim Lemma 3.6.  Let $\Omega$ satisfy the hypotheses of Lemma 3.5.
If $\eps < {{\eps_0}/ 2}$, then the
hyperbolic geodesic $\gamma$ from $0$ to $1$ for the region $\Om$ lies
in the kite
$$P=\{z: |\arg z|< \eps\}\cap\{z: |\arg(1-z)| < {{5\eps}\over 6}\},$$
and the tangent vectors to $\gamma$ have argument less than ${8\over
3}\eps$.
\bigskip

\proof By J\o rgensen's theorem, $\gamma$ is contained in the closed
disc through $1$  and $0$ which has slope $\tan \eps$
at $0$. 
Likewise $\gamma$ is contained in the reflection of this disc about
$\R$ and hence $|\arg z| < \eps$ on $\gamma$. This also shows that
$\gamma$ is contained in a kite like $P$ but with angles
$\eps$ at both $0$ and $1$. In the proof of Theorem 3.4, however,
we need the
improvement to ${{5\eps}\over 6}$ of the angle at $1$. 

By Lemma 3.5,
a portion of $\gamma$ near $1$ lies in $P$. 
Suppose $w_1\in \gamma$ with 
$|\arg w_1| =\delta < \eps$ and then apply Lemma 3.5 to the
region ${1\over{w_1}}\Om$ with $\eps$ replaced by $\eps +\delta$. Then
the tangent vector to $\gamma$ at $w_1$ has argument $\theta$ where

$$|\theta - \arg w_1| < {5\over 6}(\eps + |\arg w_1|).\eqno(3.3)$$
Since $|\arg w_1| < \eps$, we have $|\theta| \le {8\over 3}\eps.$
Moreover (3.3) also implies 
$\theta < {5\over 6}\eps$, when $\arg w_1 \le 0$, and
$\theta  > -{5\over 6} \eps$ when $\arg w_1 \ge 0$.
But if $w_1$ is the last point on $\gamma\cap \b P$ before reaching
$1$, this is impossible. Thus $\gamma \subset P$, proving the lemma.
\endproof

\bigskip

The next lemma improves Lemma 3.5 by only requiring that
the portion of $\b\Om$ 
in a large disc lies inside a small sector.

\bigskip
\proclaim Lemma 3.7.
There is a constant $C_1$ so that if $\eps < {{\eps_0}/ 2}$ and 
if $\b\Om$ is a Jordan arc  such that
$0\in \b\Om$, $\b\Om \cap \{|z|> C_1/\eps^2\}\ne\emptyset$, and
$$\{z: |\arg z| < \pi -\eps \hbox{ and } |z|\le {{C_1}\over \eps^2} \}
\subset \Om,$$ 
then  the conformal map $f:\H^+ \longrightarrow \Om$ with $f(0)=0$ and 
$|f(\infty)|> {{C_1}\over {\eps^2}}$ satisfies
$$|\arg z_0^2 f'(z_0)| < {{9\eps}\over 10},\eqno(3.4)$$
where $z_0=f^{-1}(1)$.
\bigskip

\proof As before, we may assume $|z_0|=1$.  Let $\om(z, E, V)$ denote
harmonic measure at $z$ for $E\cap \overline{V}$ in $V\setminus E$. 
Set $R= {{C_1}\over {\eps^2}}$ and $B_R=B(0,R)=\{|z|<R\}$.
Then by
Beurling's projection theorem
and a direct computation 
$$\om (1,\b B_R, \Om) \le \om(1,\b B_R ,B_R\setminus[-R,0])={4\over{\pi}}
\tan^{-1}\biggl({1\over{R^{1\over 2}}}\biggr).\eqno(3.5)$$
By the maximum principle
$$\biggl|\arg{{f(z_0)}\over{z_0^2}}\biggr|
\le \eps + (2\pi + \eps){4\over{\pi}}
\tan^{-1}\biggl({{\eps}\over{C_1^{1\over 2}}}\biggr) < {{11\eps}\over
10},$$
for $C_1$ sufficiently large. Since $f(z_0)=1$, 
$$|\arg {{z_0}}|\le {{11\eps}\over {20}}. \eqno(3.6)$$

Next we show that there is a large half disc contained in
$f^{-1}(\Om\cap B_R)$. Set 
$$S= \inf\{|w-i\Im z_0|: w\in \H^+ \hbox{ and } f(w) \in \b B_R\}.$$
Using the map 
$${{z-i \Im z_0-S}\over{z-i\Im z_0 + S}}$$
of $\H^+$ onto $\D$  and Beurling's projection theorem again,
$$\om(z_0, f^{-1}(\b B_R),\H^+) \ge \om(z_0,[S,\infty)+i \Im z_0,\H^+).$$ 
\noindent Then by (3.5), (3.6) and an explicit computation
$${{4}\over{\pi}}\tan^{-1}\biggl({{\eps}\over{C_1^{1\over 2}}}\biggr)
\ge {{2}\over{\pi}}\tan^{-1}\biggl({{\Re z_0}\over{\sqrt{S^2-\Re
z_0^2}}}\biggr). $$
For $\eps$ sufficiently small, this implies
$$B(0,{{C_1^{1\over 2}}\over{2\eps}})\cap \H^+ \subset 
f^{-1}\bigl(\Om \cap B(0,{{C_1}\over {\eps^2}})\bigr).$$
 
Now follow the proof of Lemma 3.5 replacing $\tau$ with a conformal
map of $\D$ onto $\H^+\cap \{|z|< {{C_1^{1\over 2}}\over {2\eps}}\}$
such that $\tau(0)=z_0$. Then $\tau'(0)= 2\Re z_0 +
O({{\eps}\over{C_1^{1\over 2}}})$ and for $C_1$ sufficiently large,
(3.4) holds. 
\endproof

\bigskip
Following the proof of Lemma 3.6 (replacing 5/6 by 9/10), the next corollary obtains.

\bigskip
\proclaim Corollary 3.8. 
Suppose $\b\Om$ is a Jordan arc  such that
$0\in \b\Om$, $\b\Om \cap \{|z|> C_1/\eps^2\}\ne\emptyset$, and
$$\{z: |\arg z| < \pi -\eps \hbox{ and } |z|\le {{C_1}\over \eps^2} \}
\subset \Om.$$ 
If $\eps < {{\eps_0}/ 2}$, then the
hyperbolic geodesic $\gamma$ from $0$ to $1$ for the region $\Om$ lies
in the kite
$$P=\{z: |\arg z|\le \eps\}\cap\{z: |\arg(1-z)|\le {{9\eps}\over
10}\}.\eqno(3.7)$$
Moreover, the tangent vectors to this geodesic have argument at most $3\eps$.
\bigskip 

\noindent{\bf Proof of Theorem 3.4.} 
Use the constant $C_1$ from Lemma 3.7 in Definition 3.3.
As in Theorem 2.2, let $\gamma_j$ denote the portion of the 
computed boundary $\b\Om$ between $z_j$ and $z_{j+1}$. By construction 
$\gamma_0\cup\gamma_1$ is a half line through $z_0=\infty$, $z_1$, and
$z_2$. Make the inductive hypotheses that 
$$\bigcup_{j=0}^{k-1}\gamma_j \subset \bigcup_{j=0}^{k-1}
D(z_j,z_{j+1})\eqno(3.8)$$ 
and
$$\gamma_{k-1}\cap P_k =\emptyset.\eqno(3.9)$$
Since the $\eps$-diamond chain 
$D(\infty,z_1), D(z_1,z_2), \dots D(z_{n-1},z_n)$
satisfies the $\eps$-pacman condition, (3.8) and (3.9) show that the
hypotheses of  
Corollary 3.8 hold for the curve $\gamma=\cup_0^{k-1}\gamma_j$ and hence
$\gamma_k\subset D(z_k,z_{k+1})$. 
Also by Corollary 3.8 and (3.1), 
$$\gamma_k \cap P_{k+1}=\emptyset.$$
By induction, the theorem follows.
\endproof

If the hypotheses of Theorem 3.4 hold, then the proof of Proposition
2.5 gives the following Corollary.
\bigskip

\proclaim{Corollary 3.9}. If $\Om$ and the diamond chain
$D(z_k,z_{k+1})$ satisfy the hypotheses of Theorem 3.4, then the
conformal map $\f_c$ computed in the geodesic algorithm extends to be
conformal on $\Om\cup\bigcup_{k=0}^n D(z_k,z_{k+1})$.

The next Theorem says that for a region $\Om$ bounded by a $C^1$ curve, 
the geodesic algorithm with data points $z_0, z_1,\dots, z_n$ 
produces a region $\Om_c$ whose boundary is
a $C^1$ approximation to $\b \Om$. 

\bigskip
\proclaim Theorem 3.10. Suppose ~$\Om$ is a Jordan region bounded by a $C^1$
curve ~$\b \Om$. Then there exists $\delta_0>0$, depending on $\b\Om$
so that for $\delta<\delta_0$, ~$\b \Om$ is contained in a closed
$\delta$-diamond-chain $D=\cup D(z_k,z_{k+1})$  and so that
$\b\Om_c$, the boundary of the region computed by the geodesic algorithm, 
is contained $D$. Moreover
if $\zeta\in \b \Om_c$ and if
$\alpha\in \b \Om$ with $|\zeta-\alpha| < \delta$ then
$$|\eta_\zeta - \eta_\alpha| < 6\delta,\eqno(3.10)$$
where $\eta_\zeta$ and $\eta_\alpha$ are the unit tangent vectors to
$\b \Om$ and $\b \Om_c$ at $\zeta$ and $\alpha$, respectively.
\bigskip

\proof There were two reasons for requiring that $z_0=\infty$
in Theorem 3.4. The first reason was to assure that 
$$\Bigl(\cup_0^{k-1} \gamma_j\Bigr)\cap (\C\setminus B(z_k,R_k))\ne
\emptyset\eqno(3.11)$$
as needed for Lemma 3.7. 
The second reason is the difficulty in closing the curve, since Lemma
3.7 does not apply.
The difficulty being that a pacman centered
at $z_n$ will contain $z_0$ if $z_0$ is too close to $z_n$. 
Since $\b \Om\in  C^1$, we may
suppose that the $\delta$-diamond chain $D(z_0,z_1), D(z_1,z_2),\dots,
D(z_{n-1},z_n)$ satisfies the pacman condition. Note that this
requires $z_n$ to be much closer to $z_{n-1}$ than to $z_0$.
Since $\b \Om \in C^1$, if $|z_n-z_0|$ is sufficiently small, 
we can find two discs 
$$\Delta_p\subset \C\setminus \bigcup_0^{n-1} D(z_k,z_{k+1}),$$ 
for $p=1,2$, 
with 
$$\{z_0, z_n\}= \b \Delta_1 \cap \b \Delta_2\subset
\Delta_1\cap\Delta_2\subset
D(z_n,z_0),$$
where $D(z_n,z_0)$ is a $\delta$-diamond.
By J\o rgensen's theorem, as in the proof of Theorem 2.2,
the geodesic $\gamma_n$ from $z_n$ to
$z_0$ is contained in $\Delta_1\cap \Delta_2$. 
Then by the proof of Theorem 3.4,
$\b \Om_c$ is contained
in the $\delta$-diamond chain. The statement about tangent vectors now
follows from Corollary 3.8.
\endproof

\bigskip

We say that $\{z_k\}$ are {\bf locally evenly spaced} 
if 
$${1\over D} \le \biggl|{{z_{k}-z_{k-1}}\over{z_k-z_{k+1}}}\biggl|
\le D, \eqno(3.12)$$
for some constant $D<\infty$.
Note that the spacing between points can
still grow or decay geometrically.  
We define the {\bf mesh size}
$\mu$ of
the data points $\{z_j\}$ to be
$$\mu(\{z_j\})=\sup_k |z_k-z_{k+1}|.$$

We say that a Jordan curve $\Gamma$  in the extended plane $\C^*$
is a {\bf $K$-quasicircle} if for some linear fractional transformation
$\tau$
$${{|w_1-w|+|w-w_2|}\over{|w_1-w_2|}}\le K\eqno(3.13)$$ 
for all $w_1, w_2 \in \tau(\Gamma)$ and for all $w$ on
the subarc of $\tau(\Gamma)$ with smaller diameter.
Thus circles and lines are $1$-quasicircles. Quasicircles look very
flat on all scales if $K$ is close to $1$, but for any $K>1$ they 
can contain  a 
a dense set of spirals. See for example, Figure 8.

If $\Gamma$ satisfies (3.13) with $K=1+\delta$ and small $\delta$ and if
$\{z_k\}\subset\tau(\Gamma)$ is locally evenly spaced then
$$\biggl|\arg\biggl({{z_k-z_{k-1}}\over{z_{k+1}-z_k}}\biggl)\biggl|\le
C\delta^{1\over 2},\eqno(3.14)$$
for some constant $C$, depending on $D$.


\bigskip
\proclaim Theorem 3.11. There is a constant $K_0>1$ so that if
$\Gamma$ is a $K$-quasicircle with $K=1+\delta < K_0$ and if $\{z_k\}$
are locally evenly spaced on $\Gamma$, then the geodesic
algorithm finds a conformal map of $\H$ onto a region $\Om_c$ bounded by a
$C(K)$-quasicircle containing the data points $\{z_k\}$. 
The constant
$C(K)$ can be chosen so that $C(K) \to 1$ as $K\to 1$.
Moreover, given
$\eta>0$, if the mesh size $\mu(\{z_k\})$ is sufficiently small then 
$$d_H(\Gamma,\b\Om_c) < \eta,$$
where $d_H$ is the Hausdorff distance in the spherical metric. 

\proof We may suppose that $\Gamma$ satisfies (3.13) 
with $K=1+\delta$ and $\delta$ small. Note that $\infty\in \Gamma$.
If $\{z_k\}_1^n$ are locally evenly spaced points on 
$\b \Om$, with $\mu=\max|z_k-z_{k-1}|$ sufficiently small then 
(3.14) holds and
$D(\infty,z_1),
D(z_1,z_2),...,D(z_{n-1},z_n), D(z_n,\infty)$ is an $C\delta^{1\over
2}$-diamond chain, where the main axis of  the cone $D(\infty,z_1)$ is in the
direction $z_1-z_2$ and the main axis of  $D(z_n,\infty)$ is in the
direction $z_n-z_{n-1}$.
 Moreover 
$D(\infty,z_1), D(z_1,z_2),...,D(z_{n-1},z_n)$ 
satisfies the $\eps$-pacman condition if
$$\eps \ge C \delta^{1\over 4},$$ 
for some universal constant $C$.  Now apply Theorem 3.4 to obtain
$\gamma_j\subset D(z_{j-1},z_j)$, $j=1,\dots,n-1$. 
By an argument similar to the proof of Theorem 3.10, we 
can also find a geodesic arc for $\C\setminus (\cup_0^{n-1} \gamma_j)$
 from $z_n$ to $\infty$ contained in 
$D(z_n,\infty)$. Then the computed curve will be a
$CK$-quasicircle. 
\endproof
\bigskip
As noted before, the boundary of the region computed with the geodesic
algorithm, $\b\Om_c$, is a $C^1$ curve. 
We end this
section by proving that $\b\Om_c$ is slightly better than $C^1$.
If $0<\alpha<1$, we say that a curve $\Gamma$ {\bf belongs to}
${\bf C^{1+\alpha}}$ if 
arc length parameterization 
$\gamma(s)$ of  
$\Gamma$ satisfies
$$|\gamma'(s_1)-\gamma'(s_2)|\le C|s_1-s_2|^{\alpha}$$
for some constant $C<\infty$. 

\bigskip
We say that a conformal map $f$ defined
on a region $\Om$ {\bf belongs to} ${\bf C^{1+\alpha}(\overline{\Om})}$,
~${0<\alpha <1}$
provided $f$ and $f'$ extend to be 
continuous on $\overline\Om$ and there is a constant $C$ so that
$$|f'(z_1)-f'(z_2)| \le C|z_1-z_2|^{\alpha}$$
for all $z_1, z_2$ in $\overline{\Om}$.

\bigskip
\proclaim{Proposition 3.12}. If the bounded Jordan region
$\Om_c$ is the image of the unit disc
by the geodesic algorithm, then
$$\b\Om_c \in C^{3/2},$$
and $\b\Om_c \not\in C^{1+\alpha}$ for $\alpha>1/2$, unless $\Om_c$ is
a circle or a line. Moreover $\varphi \in C^{3/2}(\overline{\Om_c})$ and
$\varphi^{-1}\in  C^{3/2}(\overline{\D})$.

\bigskip

\proof
To prove the first statement,  it is enough to show
that if $\gamma$ is an arc of a circle in $\H$
which meets $\R$ orthogonally 
at $0$ (constructed by application of one of the maps $f_a^{-1}$ as in
Figure 2), then the curve $\sigma$ which is the 
image of $[-1,1]\cup \gamma$ by  the map
$S(z)=\sqrt{z^2-d^2}$ is $C^{3\over 2}$ (and no better class) 
in a neighborhood of $S(0)=i d$. Indeed,
subsequent maps in the composition $\varphi^{-1}$ are conformal in
$\H$ and hence preserve smoothness.
For $d>0$, the function 
$$\psi(z)=
\sqrt{\Biggl({{\sqrt{z^2-c^2}\over{1+\sqrt{z^2-c^2}/b}}}\Biggl)^2-d^2}=
id + {{i}\over{2d}}(z^2-c^2) - {{i}\over {bd}}(z^2-c^2)^{3\over 2} + 
{\rm O}((z^2-c^2)^2)$$
for some choice of $b\in \R$ and $c>0$ is a conformal map of the 
upper half plane
onto a region whose complement contains the curve $\sigma$.
Clearly $\psi\in C^{3\over 2}$ near $z=\pm c$, and so by a theorem of
Kellogg (see [GM, page 62]), $\sigma\in C^{3\over 2}$. The same theorem
implies $\sigma$ is not in $C^{\alpha}$ for $\alpha > {3\over 2}$
unless $1/b=0$. This argument also shows that $\f_c\in
C^{3/2}(\overline{\Om})$. To prove $\f_c^{-1}\in C^{3\over
2}(\overline{\D})$, apply the same ideas above to the inverse maps.
Alternative, this last fact can be proved by following the proof of
Lemma II.4.4 in [GM].
\endproof
\bigskip
\bigskip
\S {\bf 4. Slits and Newton's method}
\bigskip

One complication of the slit and zipper algorithms is that the
basic maps $f_a=g_a^{-1}$  cannot be written explicitly in terms of
elementary maps, unlike the geodesic algorithm. Newton's method can be
used to find the inverse of $g_a$. 

Fix $p$, with $0 < p < 1$, and let $f(z)=(z-p)^p (z+1-p)^{1-p}$.
Then $f(\H)=\H\setminus L$ where $L$ is the line segment from $0$ to
$e^{i\pi p} p^p(1-p)^{1-p}$. 
(Note that ${1\over 2}\le |L|\le 1$).
Fix $w\in f(\H)$. We wish to solve
$$f(z)=w\eqno(4.1)$$
for $z$. 
Newton's method then takes an initial guess $z_0$ and defines
$$z_{n+1}= z_n - {{f(z_n)-w}\over{f'(z_n)}}.$$
Near $\infty$
$$f(z)=z + 1-2p + {\rm O}({1\over z})$$ so a natural first guess for
an approximation to the solution $z$ to (4.1) would then be
$$z_0=w + 2p-1.$$

The next Theorem says that this initial guess $z_0$ will work for large $w$. 

\bigskip
\proclaim Theorem 4.1. If $0 < p < 1$, set
$$f(z)=(z-p)^p(z+1-p)^{1-p},$$ 
and suppose $|w| > (1+\sqrt{5})/2$.
Then for $z_0=w+2p-1$  the n-th Newton iterate $z_n$ has relative
error
$$\biggl|{{f(z_n)-w}\over{w}}\biggr| \le {3\over 2}\biggl({1\over
{12}}\biggr)^{2^{n}}.$$
\bigskip

For example $$\biggl|{{f(z_4)-w}\over{w}}\biggr| < 10^{-17}$$
so that $z_4$ is virtually a formula for $f^{-1}(w)$. In fact, in the
slit or geodesic algorithm the points $w$ with small $|w|$ correspond
to points in the region near the corresponding vertex, so that most
points will have large modulus. In practice, most points need only one
or two iterations of Newton's method. The ``approximate zero theorem''
of Smale and Shub-Smale (see [SS]) can be also used to show
that Newton's method will
converge quadratically if $|w|$ is sufficiently large.
Since we have an explicit formula for $f$, it is not surprising that
we get a somewhat stronger result, in terms of $|w|$, in Theorem 4.1.

\bigskip
\proof
Set
$F(z)=(f(z)-w)/w$. 
We claim that 
$$|F(w+2p-1)| \le {{p(1-p)}\over 2} \eqno(4.2)$$
and if
$$|F(z)| \le {{p(1-p)}\over 2}, \eqno(4.3)$$
then 
$$|z|^2\ge p(1-p).\eqno(4.4)$$

To prove these claims, we study the auxillary function
$$H(\zeta)=\bigl(1-(1-p)\zeta\bigl)^p\bigl(1+p\zeta\bigl)^{1-p} -~ 1,$$
which has derivative
$$H'(\zeta)={{p(p-1)\zeta}\over{[1+(p-1)\zeta]^{1-p}[1+p\zeta]^{p}}}.$$
Bounding the denominator from below and integrating we obtain the
estimate
$$|H(\zeta)|\le{{p(1-p)}\over 2}{{|\zeta|^2}\over
{1-|\zeta|}}. \eqno(4.5)$$

\noindent Note first that $F(w+2p-1)=H(1/w)$. So that by (4.5)
$$|F(w+2p-1)| \le {{p(1-p)}\over
2}\biggl({2\over{1+\sqrt{5}}}\biggl)^2{{1}\over{1-{2\over{1+\sqrt{5}}}}}
= {{p(1-p)}\over 2},$$
proving (4.2).

Suppose now that (4.3) holds and $|z|^2 \le p(1-p)$. Then
$$1 - \Bigl|{{(z-p)^p(z+1-p)^{1-p}}\over w}\Bigl|\le {{p(1-p)}\over
2}\le {1\over 8}.$$
This implies 
$$\eqalign{|w| &\le 
{8\over 7}[p^{1\over 2}(1-p)^{1\over 2} + p]^p
[p^{1\over 2}(1-p)^{1\over 2} + 1-p]^{1-p}\cr\cr
&={8\over 7}[p^p(1-p)^{1-p}(1+ 2p^{1\over 2}(1-p)^{1\over 2})]^{1\over
2}\cr\cr
&\le {{8\sqrt{2}}\over 7} < {{\sqrt{5}+1}\over 2},}$$
contradicting our assumption $|w| \ge (\sqrt{5}+1)/2$, and proving
that (4.3) implies (4.4).

Next suppose that (4.3) holds and set
$$\widetilde z = z -{{F(z)}\over {F'(z)}}=z+
\biggl({{z-p}\over{z}}\biggl)\biggl[w\biggl({{z+1-p}\over{z-p}}\biggl)^p -(z+1-p)\biggl],$$
Then after some manipulations we obtain the magic formula
$$F(\widetilde z) =
H\biggl({{F(z)}\over{z}}\biggr).\eqno(4.6)$$
By (4.3) and (4.4),
$$\biggl|{{F(z)}\over z}\biggl| \le {{\sqrt{p(1-p)}}\over 2}\le
{1\over 4},$$
and so by (4.6), (4.5), (4.4) and (4.3)
$$|F(\widetilde z)| \le {{p(1-p)}\over {2({3\over
4})}}\biggl|{{F(z)}\over z}\biggl|^2 \le {2\over 3}|F(z)|^2
<{{p(1-p)}\over 2}$$
 
\noindent By induction and (4.2)
$${2\over 3}|F(z_n)| \le  \biggl({2\over 3}|F(z_0)|\biggl)^{2^n}\le
\biggl({1\over{12}}\biggl)^{2^n},$$
proving Theorem 4.1.
\endproof

%
%
%
The region of possible $w$ where quadratic convergence is obtained can
be enlarged with more involved estimates, but
Newton's method applied directly to $f$ with this initial value will
not always converge. Indeed,  the Newton interate $z-F(z)/F'(z)$ has
repelling fixed points at $p$ and $p-1$, and a pole at $0$.
In order to successfully apply the algorithm to a
wide variety of curves we need to find a reliable routine for finding
the inverse. 

In the implementation of the slit and zipper algorithms we consider
four regions based on the length $|L|= p^p(1-p)^{1-p}$ 
of the segment $L$ and the imaginary part of the tip $w_{\rm tip}= f(0)$.
\bigskip

$$\Om_\infty=\{w: |w|\ge {9\over 8} |L|\} ~~~~~~~~~~~~
\Om_{\rm tip}=\{w: |w-w_{\rm tip}| < {1\over 4} \Im w_{\rm tip}\}$$

$$\Om_p =\{w: 0 < \arg w < \pi p\} ~~~~~~~~~~~~
\Om_{p-1} = \{w: \pi p < \arg w < \pi\}$$
\bigskip

\vskip 0.5truein
\psfrag{Omi}{$\Om_{\infty}$}
\psfrag{Oma}{$\Om_p$}
\psfrag{Omoa}{$\Om_{p-1}$}
\psfrag{Omt}{$\Om_{tip}$}
\psfrag{pa}{$\pi p$}
\psfrag{R}{$\R$}
\psfrag{L}{$L$}
\centerline{\includegraphics[height=1.5in]{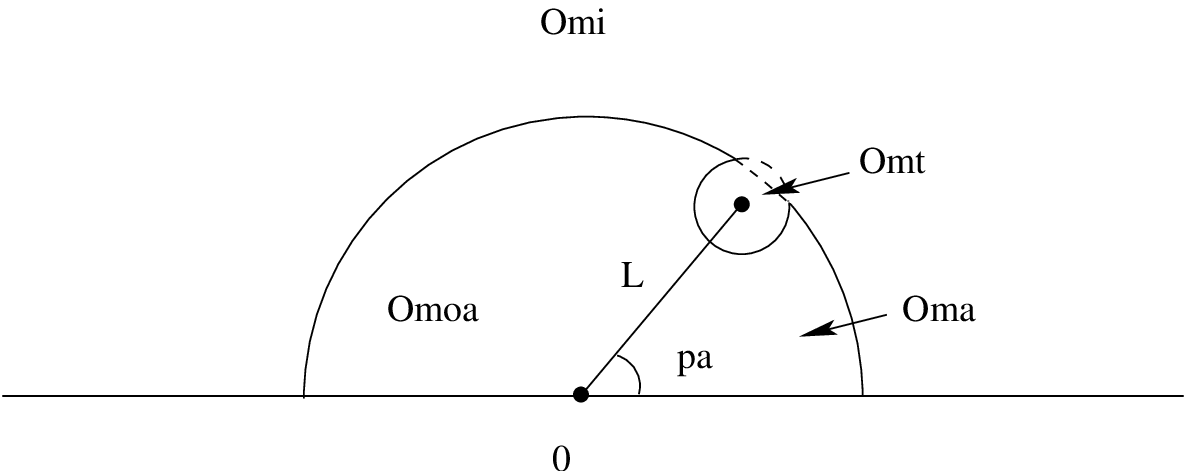}}
\nobreak
\centerline{{\bf Figure 10.} ~Iteration regions for Newton's method}
\vskip 0.5truein

If $w\in \Om_{\infty}$ then we use Newton's method as described in
Theorem 4.1. However, we improve the initial guess by taking more
terms in the expansion at $\infty$:
$$z_0= w + {{2p-1}} + {{p(1-p)}\over{2w}} +
{{(1-2p)p(1-p)}\over{3w^2}},$$
and we rewrite the function to iterate on $z/w$ instead of $z$ to
improve numerical accuracy.
If $w\notin \Om_{\infty}$ but $w\in \Om_{\rm tip}$ then
we first open up the region by applying $k(w)=\sqrt{w-w_{\rm tip}}$.
Then $k\circ f$ extends to be analytic and one-to-one in a
neighborhood of $0$. So we use Newton's method to solve $k\circ
f(z)= k(w)$ for $z$. The remaining $w$ are in the sectors between
$\R$ and $L$. If $w \notin \Om_\infty\cup \Om_{\rm tip}$
but $0 < \arg w < \pi p$, then we apply the preliminary map $k_p(z) =
z^{1 \over p}$ instead of $k$ and use Newton's method again. For the
remaining points we use the preliminary map $k_{1-p}(z)
=z^{1\over{1-p}}.$  
We leave the proof of the analog of Theorem 4.1 in
the remaining three regions to the interested reader. 
While we cannot prove convergence of Newton's
method in every case, extensive numerical testing
indicates that we have chosen the proper regions.

\bigskip
\bigskip
\S {\bf 5. Estimates for conformal maps onto nearby domains}
\bigskip

We begin this section with a discussion of the following question. 
Consider two simply connected planar domains
$\Om_j$ with $0\in \Om_j$ and conformal maps $\varphi_j:\Om_j\to\D$ fixing 0,
suitably normalized (for instance positive derivative at 0). 
If $\Om_1$ and $\Om_2$ are ``close,'' what can be said about $|\f_1-\f_2|$ on $\Om_1\cap \Om_2$, 
or about $|\f^{-1}_1-\f^{-1}_2|$ on $\D$?
The article [W] gives an overview and numerous results
in this direction. How should ``closeness'' of the two domains be measured?
Simple examples show that the Hausdorff distance in the Euclidean or spherical metric
between the boundaries does not give uniform estimates for either $||\f_1-\f_2||_\infty$ or
$||\f_1^{-1}-\f_2^{-1}||_\infty$. 

\vskip 0.5truein
\psfrag{Oa}{$\Om_1$}
\psfrag{z1}{$z_1$}
\psfrag{z2}{$z_2$}
\psfrag{z3}{$z_3$}
\centerline{\includegraphics[height=1.5in]{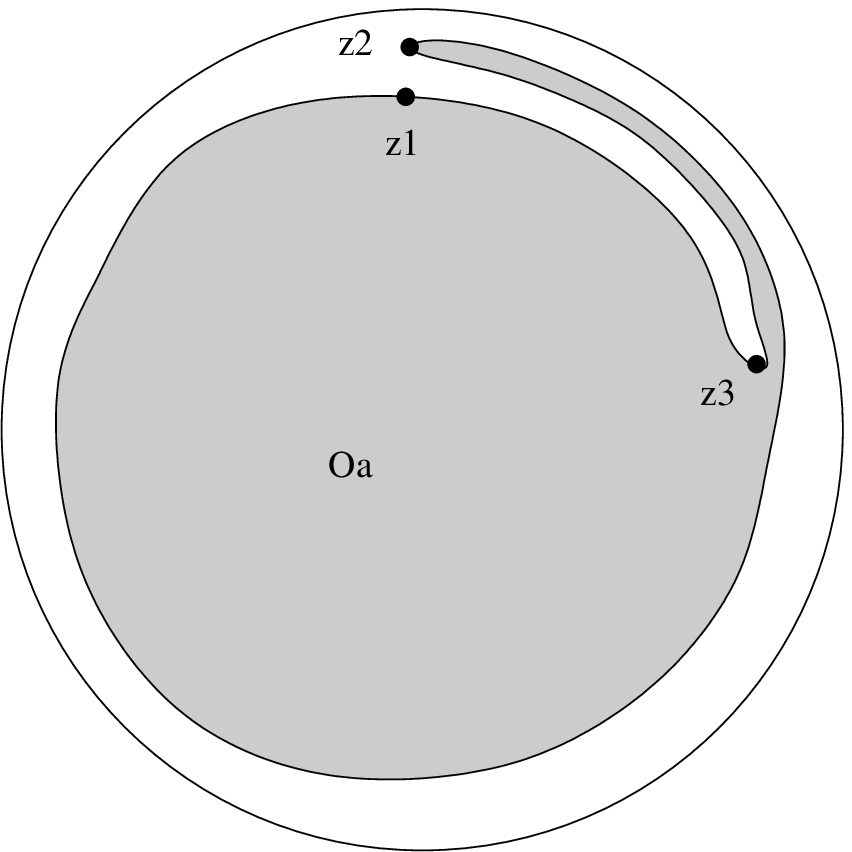}}
\nobreak
\centerline{{\bf Figure 11.} ~Small Hausdorff distance }
\vskip 0.5truein

\noindent For example in Figure 11, $\Om_1$ contains a disc of radius $1-\delta$ where $\delta$ is small
and hence 
for $\Om_2=\D$, ~$d_H(\Om_1,\Om_2)\le \delta$, but 
$|\f_1(z_1)-\f_1(z_2)|$ is large and $|\f_1(z_2)-\f_1(z_3)|$ is small so
that neither $||\f_1(z)-z||_\infty$ nor $||\f_1^{-1}(z)-z||_\infty$ is small.

Mainly for ease of notation, we will assume throughout this section that the $\Om_j$ are Jordan domains,
and denote $\gamma_j:\bD \to\partial \Om_j$ an orientation preserving parametrization.
Even the more refined distance
$$\inf_{\alpha} ||\gamma_1 -\gamma_2\circ\alpha||_{\infty},$$
where the infimum is over all homeomorphisms 
$\alpha$ of $\bD,$ does not control $||\f^{-1}_1-\f^{-1}_2||_{\infty}$
or $||\f_1-\f_2||_{\infty}$. For example, let $\Om_2$ be a small rotation of the region $\Om_1$ in
Figure 11. What is needed is some control on the ``roughness'' of the boundary.
Following [W], for a simply connected domain $\Om$ we define
$$\eta(\delta) = \eta_\Om(\delta) = \sup_{C} \diam \ T(C),$$
where the supremum
 is over all crosscuts of $\Om$ with $\diam \ C\leq \delta$, and where
$T(C)$ is the component of $\Om\setminus C$ that does not contain 0.
Notice that $\eta(\delta)\to 0$ as $\delta\to 0$ is equivalent to saying that $\partial \Om$ is locally
connected, and the condition $\eta(\delta)\leq K\delta$ for some constant $K$ is equivalent
to saying that $\Om$ is a John-domain (e.g. [P], Chapter 5).
It is not difficult to control the modulus of continuity of $\f^{-1}:\D\to \Om$ in term of $\eta,$
see [W], Theorem I. This can be used to estimate $||\f^{-1}_1-\f^{-1}_2||_{\infty}$ in terms of
the Hausdorff distance between the boundaries, for example.

\bigskip
\proclaim Theorem 5.1 (Warschawski[W], Theorem VI). If $\Om_1$ and $\Om_2$ are John-domains, 
$\eta_j(\delta)\leq \kappa \delta$ for $j=1,2$, and if $d_H(\partial \Om_1, \partial \Om_2)\leq \epsilon$, 
then
$$||\f^{-1}_1-\f^{-1}_2||_{\infty}\leq C \epsilon^{\alpha}$$
with $\alpha=\alpha(\kappa)$ and $C=C(\kappa,\dist(0,\partial \Om_1\cup \partial \Om_2))$.
\bigskip

In fact, Warschawski proves that every $\alpha < 2/(\pi^2 \kappa^2)$ will work (with $C=C(\alpha)$). 
Using the H\"older continuity of quasiconformal maps, his proof can easily be modified to give the 
following better estimate 
if $\Om_1$ and $\Om_2$ are $K$-quasidiscs with $K$ near 1. A
$K$-quasidisc is a Jordan region bounded by a $K$-quasicircle.

\bigskip
\proclaim Corollary 5.2.
If $\Om_1$ and $\Om_2$ are $K$-quasidiscs, and if $d_H(\partial \Om_1, \partial \Om_2)\leq \epsilon$, 
then
$$||\f^{-1}_1-\f^{-1}_2||_{\infty}\leq C \epsilon^{\alpha}$$
with $\alpha=\alpha(K)\to 1$ as $K\to1.$
\bigskip

As for estimates of $||\f_1-\f_2||_{\infty}$, Warschawski shows [W, Theorem VII] that

$$\sup_{\Om_1}|\f_1-\f_2|\leq C \epsilon^{1/2}\log{2\over{\epsilon}}$$
if $\Om_1\subset \Om_2$, and if $\Om_1$ is a John-domain, with $C$ depending on $\kappa$
and on $\dist(0,\partial \Om_1\cup \partial \Om_2)$.
However, his result does not apply without the assumption of inclusion $\Om_1\subset \Om_2$.
To treat the general case 
the trick of controlling $|\f_1-\f_2|$ by passing to the conformal map $\f$ of the component
$\Om$ of $\Om_1\cap \Om_2$ containing $0$ (which now is included in $\Om_j$) 
does not seem to work, as the geometry
of $\Om$ can not be controlled.
Nevertheless, for the case of disc-chain domains, the above estimate can be proved, even
without any further assumption on the geometry on the circle chain:

\bigskip
\proclaim Theorem 5.3.
Let $D_1, D_2,...,D_n$ be a closed $\eps$-disc-chain 
surrounding 0.
Suppose $\partial \Om_j\subset \cup_k \overline{D_k}$ for $j=1,2,$ and let $\f_j:\Om_j \to \D$ be conformal
maps with  $\f_1(0)=\f_2(0)=0$ and $\f_1(p)=\f_2(p)$ for
a point $p\in \partial \Om_1\cap \partial \Om_2$. Then
$$\sup_{w\in\Om_1\cap\Om_2}|\f_1(w)-\f_2(w)|\leq C \epsilon^{1/2}\log{1\over{\epsilon}},$$
where $C$ depends on $\dist(0, \cup_k D_k)$ only.
\bigskip

In case we have control on the geometry of the domains, we have the following counterpart to
Corollary 5.2.

\bigskip
\proclaim Theorem 5.4.
If $\Om_1$ and $\Om_2$ are $K$-quasidiscs, if 
$d_H(\partial \Om_1, \partial \Om_2)\leq \epsilon$, 
and if $\f_1(p_1)=\f_2(p_2)$ for a pair of points $p_j\in\partial \Om_j$
with $|p_1-p_2|\leq\epsilon$,
then
$$\sup_{w\in\Om}|\f_1(w)-\f_2(w)|\leq C \epsilon^{\alpha}$$
with $\alpha=\alpha(K)\to 1$ as $K\to1$, where $\Om$ is the component
of $\Om_1\cap\Om_2$ containing $0$.
\bigskip

The proofs of both theorems rely on the following harmonic measure estimate,
which is an immediate consequence of a theorem of Marchenko [M]
(see [W, Section 3], for the statement and a proof). To keep this paper self-contained,
we include a simple proof, shown to us by John Garnett, for which we thank him.

\bigskip
\proclaim Lemma 5.5.
Let $0<\theta<\pi$, $0<\epsilon<1/2$ and set
$D=\D\setminus\{r e^{i t}: -\theta \leq t \leq \theta, ~~1-\epsilon \leq r < 1\}$,
$A=\partial D \setminus \bD.$
Then 
$$\omega(0,A,D)\leq {{\theta}\over{\pi}} + C \epsilon \log {1\over{\epsilon}}$$
for some universal constant $C.$
\bigskip

\proof Set $\omega(z)=\omega(z,A,D)$ for $z\in D.$ By the mean value property, it is enough to 
show that
$$\omega(z) \leq C {{\epsilon}\over{t-\theta}}$$
for $z=(1-\epsilon) e^{i t}$ and $\theta+\epsilon \leq t\leq \pi.$
To this end, set $I= \{ e^{i \tau}: -\theta \leq \tau \leq \theta \}$
and consider the circular arc $\{\zeta : \omega(\zeta,I,\D)= {1\over 3}\}$.
If $\epsilon < \epsilon_0$ for some universal $\epsilon_0$
(for $\epsilon \geq \epsilon_0$ there is nothing to prove), then
$A$ is disjoint from this arc and it follows 
that $\omega(\zeta,I,\D)\geq {1\over 3}$ on $A$.
The maximum principle implies 
$\omega(\zeta)\leq 3 \omega(\zeta,I,\D)$ on $D$. Now the desired inequality follows from
$$\omega((1-\epsilon) e^{i t},I,\D) = {1\over{2\pi}}\int_{-\theta}^{\theta} 
{1-(1-\epsilon)^2 \over |(1-\epsilon) e^{i t}-e^{i\tau}|^2} d\tau
\leq C\epsilon \int_{-\theta}^{\theta} {1\over(t-\tau)^2} d\tau <
C {{\epsilon}\over{t-\theta}}.\eqno{\endpf}$$

\bigskip
\noindent{\bf Proof of Theorem 5.3.} We may assume that $\f_j(p)=1.$
We will first assume that $p$ is one of the points $D_k\cap D_{k+1}$.
Denote $\Om$ the largest simply connected domain $\subset\C$ containing $0$ whose boundary is contained in 
$\cup_k D_k$
(thus $\overline\Om$ is the union of $\cup_k D_k$ and the bounded component of $\C\setminus \cup_k D_k$), 
and $\f$ the conformal map from $\Om$ to $\D$ with $\f(0)=0$ and $\f(p)=1.$
First, let $z\in \partial \Om_1 \cap \partial \Om$.
Denote $B$ respectively $B_1$  the arc of $\partial \Om$ ($\partial \Om_1$) from $p$ to $z.$
By the Beurling projection theorem (or the distortion theorem), 
every $\f(D_j)$ has diameter $\leq C \sqrt{\epsilon}$.
Therefore $\f(B_1)$ is an arc in $\overline\D$, with same endpoints as $\f(B)$, that is
contained in 
$S=\{r e^{it}: 1-C \sqrt{\epsilon} \leq r < 1, -C \sqrt{\epsilon} < t < \arg \f(z) + C \sqrt{\epsilon} \}.$
Denote $A=\partial S.$ 
By Lemma 5.5, 
$$\omega(0,B_1,\Om_1) \leq \omega(0,B_1, \Om\setminus B_1) \leq 
\omega(0,A,\D\setminus A)
\leq {1\over{2\pi}}\arg \f(z) + 2 C \sqrt{\epsilon} + C \sqrt{\epsilon} \log {1\over \sqrt{\epsilon}}$$
and we obtain 
$$\arg \f_1(z)=2\pi\om(0,B_1,\Om_1) \leq \arg \f(z) + C \epsilon^{1/2}\log{1\over{\epsilon}}.$$
The same argument, applied to the other arc from $p$ to $z$, gives the opposite inequality, and together it follows that
$$|\f(z)-\f_1(z)|\leq C \epsilon^{1/2}\log{1\over{\epsilon}}.$$

Now let $z\in\partial \Om_1$ be arbitrary.
If $z'$ is a point of $\partial \Om_1 \cap \partial \Om$ in the same disc $D_j$ as $z$, then we have
$$|\f(z)-\f_1(z)|\leq |\f(z)-\f(z')| + |\f(z')-\f_1(z')| + |\f_1(z)-\f_1(z')|\leq 2 C \sqrt{\epsilon}+ C \epsilon^{1/2}\log{1\over{\epsilon}}.$$ 
The maximum principle yields 
$|\f-\f_1|\leq C \epsilon^{1/2}\log{1\over{\epsilon}}$ on
$\Om_1$. The
same argument applies to $|\f-\f_2|$,
and the theorem follows from the triangle inequality.

If $p\in \partial \Om_1\cap \partial \Om_2$ is arbitrary, let $p'$ be one of the points
$D_k\cap D_{k+1}$ in the same disc $D_j$ as $p$. Then the above estimate, applied to a rotation of $\f_1, \f_2$
and $p'$ gives
$|\f_2(p')/\f_1(p') \f_1-\f_2|\leq C \epsilon^{1/2}\log{2\over{\epsilon}}$
and the theorem follows from $|\f_j(p)-\f_j(p')|\leq C \sqrt{\epsilon}.$
\endproof

\bigskip
\noindent
The following lemma is another easy consequence of the aforementioned theorem of Marchenko [M]
([W], Section 3).

\bigskip
\proclaim Lemma 5.6.
Let $H\subset\D$ be a $K$-quasidisc with $0\in H$ such that $\partial H\subset \{1-\epsilon < |z| < 1\},$
and let $h$ be a conformal map from $\D$ to $H$ with $h(0)=0$ and $|h(p)-p|<\epsilon$ for some
$p\in\bD.$ 
Then
$$|h(z)-z| \leq C \epsilon \log {1\over{\epsilon}},$$
where $C$ depends on $K$ only.
\bigskip

\proof  We may assume $p=1.$
Let $z=e^{i \tau}$ and consider the arc $A=\{h(e^{i t}):0\leq t\leq \tau\} \subset \partial H$ of 
harmonic measure $\tau/2\pi.$ For suitable $C=C(K)$ we have that
$D=\D\setminus\{r e^{i t}: - C\epsilon \leq t \leq \arg h(z) + C\epsilon, 1-\epsilon \leq r < 1\}$
contains $A$. By the maximum principle and Lemma 5.5,
$$\tau/2\pi = \omega(0,A,H) \leq \omega(0,\partial D\cap\D,D)\leq \arg h(z)/2\pi + C \epsilon \log {1\over{\epsilon}}.$$
Applying the same reasoning to $\partial H \setminus A,$ the lemma follows for all $z\in\bD$ and thus for all $z\in\D.$
\endproof
Note that the conclusion of Lemma 5.6 is true if instead of assuming
$H$ is a $K$-quasidisc, we only assume $\arg z$ is increasing on $\b
H$.
\bigskip
\noindent{\bf Proof of Theorem 5.4.}
Because $\Om_1$ and $\Om_2$ are $K$-quasidiscs, $\f_1$ and $\f_2$ have $K^2$-quasiconformal extensions to $\C$ (see [L], Chapter I.6). 
In particular,
they are H\"older continuous with exponent $1/K^2$ (see [A]), 
and it follows that with $\alpha=1/K^2$ and $r=1-C \epsilon^{\alpha},$
we have $\f^{-1}_1(\{|z|\leq r\}) \subset \Om_2$. In particular,
$h(z) = \f_2(\f^{-1}_1(r z))$ is a conformal map from $\D$ onto a $K^4$-quasidisc $H\subset\D$, and by the H\"older continuity of
$\f_2$ and $\f^{-1}_1$ we have
$\partial H \subset \{1-C\epsilon^{\alpha^3}< |z| <1\}$. Now Lemma 5.6 yields
$|h(z)-z|\leq C\epsilon^{\beta}$, for any $\beta<\alpha^3$ and $C=C(\beta).$
For $w\in \Om\subset\Om_1\cap \Om_2,$ let $z=\f_1(w)$, then
$$|\f_1(w)-\f_2(w)|= |z-\f_2(\f^{-1}_1(z))| \leq |z-\f_2(\f^{-1}_1(r z))| + |\f_2(\f^{-1}_1(r z))-\f_2(\f^{-1}_1(z))| \leq C\epsilon^{\beta},$$
where again we have used the H\"older continuity of $\f_2$ and $\f^{-1}_1$. The Theorem follows.
\endproof

\bigskip
\bigskip
\S {\bf 6. Convergence of the Mapping Functions} 
\bigskip
We will now combine the results of Sections 2 and 3 with the estimates of the previous section,
to obtain quantitative estimates on the convergence of the geodesic algorithm.
Throughout this section, $\Omega$ will denote a given simply connected 
domain containing $0$, bounded by a Jordan curve
$\b\Om$, $z_0,..., z_n$ are consecutive points on $\b\Om$,
$\Omega_c$ is the domain and $\f_c:\Omega_c\to\D$ the map computed by the
geodesic algorithm, and $\f:\Omega\to\D$ is a conformal map, 
normalized so that $\f_c(0)=\f(0)=0$ and $\f_c(p_0)=\f(p_0)$ for some
$p_0\in \b\Om\cap\b\Om_c$.

Combining Theorems 2.2 and 5.3 and Propositions 2.5 and 3.12
we obtain at once:

\proclaim Theorem 6.1. If $\b\Om$ is  contained in a closed
$\epsilon$-disc-chain $\bigcup_{j=0}^n\overline {D_j}$
and if $z_j= \partial D_j \cap \partial D_{j+1},$
then $\b\Om_c$ is a smooth ($C^{3\over 2}$) piecewise analytic 
Jordan curve contained in $\bigcup_{j=0}^n D_j \cup z_j$, the map
$\f_c$ extends to be conformal on $\Om\cup\Om_c$ and
$$\sup_{w\in \Om} |\f(w)-\f_c(w)| \leq C \epsilon^{1/2}\log {1\over \epsilon}.$$
\bigskip

Now assume that $\b\Om$ is a $K$-quasicircle with $K<K_0$ and assume
approximate equal spacing of the $z_j$, say,
${1\over2}\epsilon < |z_{j+1}-z_j| < 2\epsilon.$ Then 
$${C\over\epsilon} \leq n \leq {C\over\epsilon^d}\eqno(6.1)$$
where $d$ (essentially the Minkowski-dimension) is close to 1 when 
$K$ is close to 1. 
Combining Theorem 3.11
with Corollary 5.2 and Theorem 5.4, we have:

\proclaim Theorem 6.2. Suppose $\b\Om$ is a $K$-quasicircle with
$K<K_0$. The Hausdorff distance between $\b\Om$ and
$\b\Om_c$ is
bounded by $C(K) \epsilon$, where $C(K)$ tends to 0 as $K$ tends to 1 and $n$ to infinity.
Furthermore,
$$||\f^{-1}-\f_c^{-1}||_{\infty}\leq C \epsilon^{\alpha}$$
and
$$\sup_{w\in \Om_0}|\f(w)-\f_c(w)|\leq C \epsilon^{\alpha}$$
with $\alpha=\alpha(K)\to 1$ as $K\to1$, where $\Om_0$ is the
component of $\Om\cap\Om_c$ containing $0$.

\bigskip
The best possible exponent in (6.1)
in terms of the standard definition of $K(\b\Om)$, which slightly differs from our geometric definition,
is given by Smirnov's (unpublished) proof of Astala's conjecture, 
$$d\leq 1+({{K-1}\over{K+1}})^2.$$
This allows us to 
easily convert estimates given in terms of $\eps$, as in Theorem 6.2,
into estimates involving $n$.

\bigskip
Finally, assume that $\b\Om$ is a smooth closed Jordan curve.  
Then $\Om$ is a $K$-quasicircle
and a John domain by the uniform continuity of the derivative of the arc length
parameterization of $\b\Om$. The quasiconformal norm $K(\b\Om)$ and
the John constant depend on the global geometry, as does the
$\eps$-pacman condition when there are not very many data points.
As the example in Figure 11 shows,  even an infinitely differentiable
boundary can have a large quasiconformal constant and a large John
constant.
However, the $\eps$-pacman condition becomes a local condition if
the mesh size $\mu(\{z_k\})=\max_k |z_{k+1}-z_k|$  of the data points
is sufficiently small.
The radii of the balls in the definition of
the $\eps$-pacman condition
$$R_k=C_1{{|z_{k+1}-z_k|}\over{\eps^2}}\eqno(6.2)$$
increase as $\eps$ decreases, but can be chosen small for a fixed
$\eps$ if the mesh size $\mu$ is small.
To apply the geodesic algorithm we suppose
that the data points have small mesh
size and, as in the
proof of Theorem 3.10, $|(z_0-z_n)/(z_{n-1}-z_n)|$ is sufficiently
large so that the $\eps$ diamond chain
$D(z_0,z_1),\dots,D(z_{n-1},z_n)$ satisfies the $\eps$-pacman
condition and 
$$\b\Om\subset\bigcup_{k=0}^n D(z_k,z_{k+1})$$
where $D(z_n,z_{n+1})=D(z_n,z_0)$ is an $\eps$-diamond.  This can be
accomplished for smooth curves by taking data points
$z_0,\dots,z_n,z_0$ with small mesh size and discarding the last few
$z_{n-n_1},\dots,z_n$ where $n_1$ is an integer depending on $\eps$
and on $\b\Om$. The remaining subset still has small mesh size (albeit
larger).
This process  of removing the last few data points 
is necessary to apply the proof of
Theorem 3.10, but in practice it is omitted. We view it only as a
defect in the method of proof.

\bigskip
If $\b\Om \in C^1$ and if $\varphi$ is a conformal map of $\Om$ onto
$\D$ then $\arg (\varphi^{-1})'$ is continuous. Indeed, it gives the
direction of the unit tangent vector. However there are examples of
$C^1$ boundaries
where $\varphi'$ and $(\varphi^{-1})'$ are not in continuous. 
In fact it is possible for both to be unbounded. If we make
the
slightly stronger assumption that $\b\Om \in C^{1+\alpha}$ for some
$0 < \alpha <1$, then $\varphi\in C^{1+\alpha}$ and 
$\varphi^{-1}\in C^{1+\alpha}$
by Kellogg's theorem (see [GM, page 62]).
In particular the derivatives are
bounded above and below  on $\overline\Om$ and $\overline\D$,
respectively.
Because of Proposition 3.12, we will consider the case $1+\alpha=3/2$.
Similar results are true for $1+\alpha < 3/2$.

\bigskip

\proclaim Theorem 6.3. Suppose $\b\Om$ is a closed Jordan curve in
$C^{3/2}$ and $\f$ is a conformal map of $\Om$ onto $\D$.  Suppose
$z_0,z_1,...z_n,z_0$
are data points on $\b\Om$ with mesh size $\mu=\max |z_j-z_{j+1}|$.
Then there is a constant $C_1$ depending on the geometry of $\b\Om$,
so that the Hausdorff distance between $\b\Om$ and $\b\Om_c$ satisfies
$$d_H(\b\Om,\b\Om_c)\le{C_1 \mu^{3/2}}\eqno(6.3)$$
and the conformal map $\f_c$ satisfies
$$||\f^{-1}-\f_c^{-1}||_{\infty}\leq   {C\mu^p}\eqno(6.4)$$
and
$$\sup_{z\in\Om\cap\Om_c} |\f(z)-\f_c(z)| \leq  {C \mu^{p}},\eqno(6.5)$$
for every $p<3/2$.

\bigskip
For example if $n$ data points are approximately evenly spaced on
$\b\Om$, so that $\mu=C/n$ then the error estimates are of the form
$C/n^{3/2}$ in (6.3) and $C/n^p$ for $p< 3/2$ in (6.4) and (6.5). 
While Theorem 6.3 gives simple estimates in terms of the mesh size or
or the number of data
points, smaller error estimates  can be obtained with fewer data points if the
data points are distributed so that there are fewer on subarcs where
$\b\Om$ is flat and more where the boundary bends or where it folds
back on itself.  In other words, construct diamond chains with angles
$\eps_k$ satisfying the $\eps_k$-pacman condition centered at $z_k$
for each $k$. The errors will then be given by 
$$\max_k \biggl(\eps_k|z_k-z_{k+1}|\biggl)^p.$$
\bigskip
\proof   
It is not hard to see from (6.2) that $\b\Om$
satisfies
the $\epsilon$-pacman condition with
$$\epsilon = {C\mu^{1/2}},$$ for $C$ sufficiently large. 
By the proof of Theorem 3.10, $\b\Om_c$ is contained in the union of the
diamonds. The diamonds $D(z_k,z_{k+1})$ have angle $C\mu^{{1/2}}$ and width
bounded by $C\mu$ and therefore (6.3) holds. 

Let $\psi$ be a conformal map of $\D$ 
onto the complement of $\overline \Om$, $\C^*\setminus
\overline{\Om}$. Then by Kellogg's Theorem as mentioned above,
$\psi\in C^{3/2}$. In particular, $|\psi'|$ is bounded above and below
on $1/2<|z|<1$. By the Koebe distortion theorem
there are constants $C_1, C_2$ so that
$$C_1(1-|z|)\le \dist(\psi(z), \b\Om)\le C_2(1-|z|),$$
for all $z$ with $1/2<|z|<1$. Thus we can choose $r= 1- C_3\mu^{3/2}$
so that the image of the circle of radius r, $I_r=\psi(\{|z|=r\})$, does
not intersect the diamond chain and 
$d_H(I_r,\b\Om) \sim \mu^{3/2}.$
Then the bounded component of the
complement of $I_r$ is a Jordan region $U_r$ containing $\Om$ and bounded by
$I_r\in C^{3/2}$, with $C^{3/2}$ norm dependent only on $\b\Om$, and
the bounds on $|\psi'|$.

Let $\sigma$ be a conformal map of $U_r$ onto $\D$. 
Inequality (6.4) now follows from [W, Theorem VIII]
by comparing the conformal
maps $\f^{-1}$ and $\f_c^{-1}$ to the conformal map $\sigma^{-1}$
where $\sigma:U_r \to \D$ and 
where all three (inverse) conformal maps are normalized to have positive
derivative at $0$ and map $0$ to the same point in $\Om$.

To see (6.5), note that

$$\sigma(\b\Om\cup\b\Om_c)\subset \{z: 1-|z| < c\mu^{3/2}\}.$$
Moreover, because $\b\Om\cup\b\Om_c$ is contained in the diamond
chain, and because both $\sigma\in C^{3/2}$ and $\sigma^{-1}\in C^{3/2}$,
$\arg \sigma(\zeta)$ is increasing along $\b\Om$,
for $\mu$ sufficiently small. By the remark after the 
proof of Lemma 5.6, 
$$|\om(0,\gamma,\sigma(\Om))-\om(0,\gamma^*,\D)| \le C \mu^{3/2} \log
\mu$$ 
for every subarc $\gamma$ of $\sigma(\b\Om)$, where $\gamma^*$ denotes the
radial projection of $\gamma$ onto $\b\D$. The same statements are
true for $\b\Om_c$. Then (6.5) follows because the harmonic measure of 
the subarc $\gamma_p$ of $\b\Om$ from $p_0$ to $p$ is given by
$$\om(0,\gamma_p,\Om)={1\over{2\pi}}\arg\biggl({{\f(p)}\over{\f(p_0)}}\biggl),$$
and a similar statement is true for $\f_c$.
\endproof

\bigskip
The constant $C$ in Theorem 6.3 depends on the
quasiconformality constant $K(\b\Om)$, $p$, $\diam(\Om)$, 
$\dist(0,\b\Om)$, and on 
$$M=\sup_{1/2<|z|<1}{\bigl(|\psi'|, 1/|\psi'|\bigl)},$$
where $\psi$ is a conformal map of the complement of $\Om$ to $\D$.
If $I_r=\psi(\{|z|=r\})$ is replaced by a $C^{3/2}$ curve
which is constructed geometrically instead of using the conformal map
$\psi$, then the constant $C$ can be taken to 
depend only on the geometry of the
region $\Om$.

\bigskip
\bigskip
\S {\bf 7. Some Numerical Results} 
\bigskip
An in depth comparision of the algorithms in this article with other methods of
conformal mapping and convergence rates 
will be written separately. To give the reader a sense of the
speed and accuracy of computations, if 10,000 data points are given, 
it takes about 20 seconds with the geodesic algorithm to compute the mapping
functions on an 3.2 GHz Pentium IV computer. Since all of the basic
maps are given explicitly in terms of elementary maps, the speed
depends only on the number of points, not the shape of the region or
the distribution of the data points.
The accuracy can be measured if the
true conformal map is known. For example 
$$f(z)= {{rz}\over {1 + (rz)^2}},$$
where $r<1$
maps the unit disc into an inverted ellipse. See Figure 12.  The
region was chosen because it almost pinches off at $0$, and because
the stretching/compression given by $\max |f'|/\min|f'|$ is 
big for $r$ near $1$.
Higher resolution images can be obtained from:
\smallskip
http://www.math.washington.edu/$\sim$marshall/preprints/zipper.pdf

\bigskip
\centerline{\includegraphics[height=1.5in]{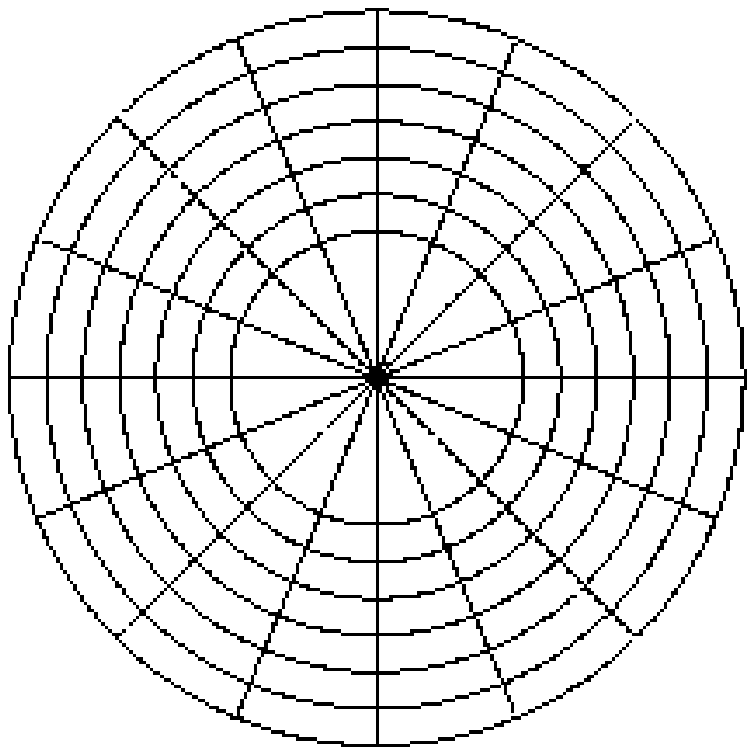}
\hfil\includegraphics[height=1.5in]{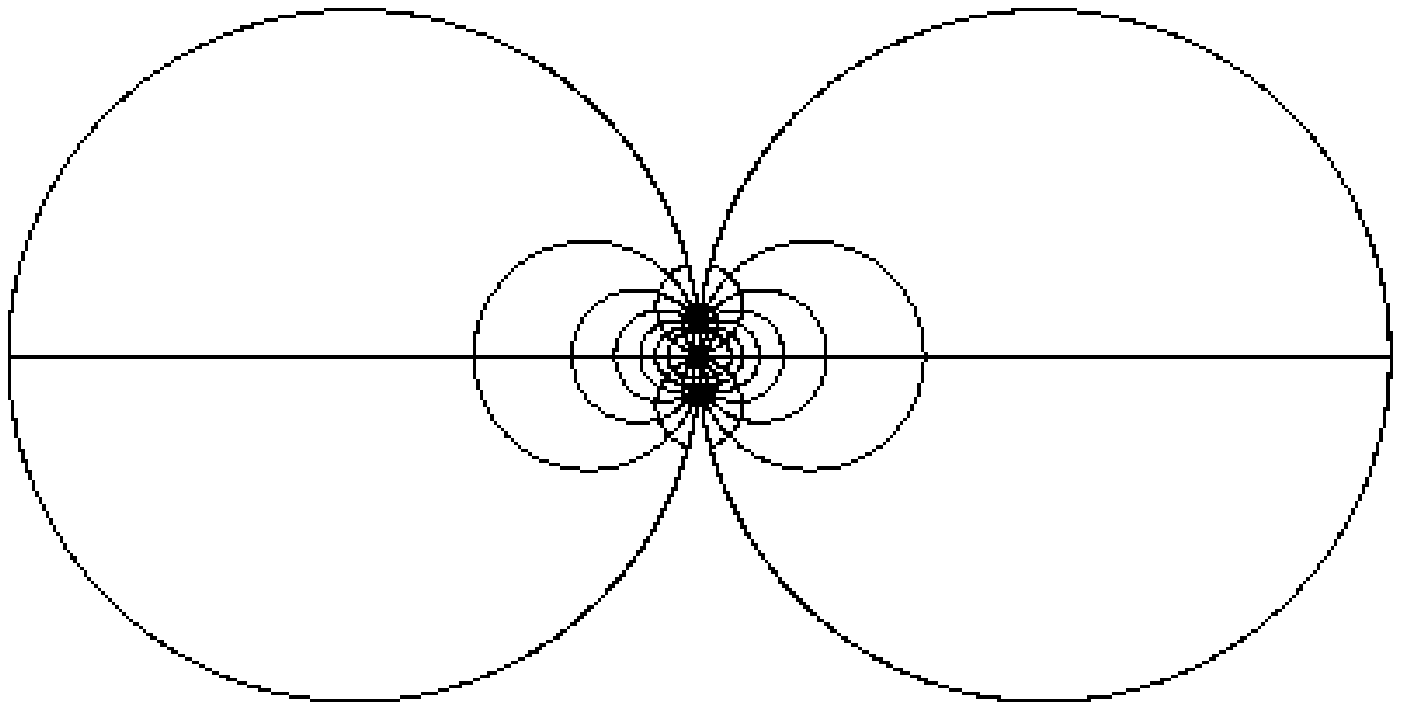}}
\nobreak
\centerline{{\bf Figure 12.} Inverted ellipse with $r=.95$.}
\bigskip

\noindent We chose $r=.95$ and
used as data points the image by $f$ of 10,000 equally spaced points 
on the unit circle, and compared the corresponding points on the unit
circle computed by the geodesic algorithm with 10,000 equally spaced
points. The errors were less than $1.8 \cdot 10^{-6}$. 
The same procedure
using the zipper algorithm took 84 seconds, and had errors less than
$9.2 \cdot 10^{-8}$.
When the number
of data points was increased to 100,000, the time to run the geodesic
algorithm increased to 25 
minutes with errors less than $2 \cdot 10^{-8}$.  In this example, the
the difference between successive boundary data points ranged from $.025$ to
$3\cdot 10^{-6}$ so that perhaps a better distribution of data points
would have given smaller errors. 

Figure 13 shows the conformal map of a Carleson grid on the disc to 
both the interior and exterior of the island Tenerife (Canary Islands).
The center of the interior is the volcano Teide.
It also shows both the  original data for the coastline, connected 
with straight line segments, and the boundary curve connecting the
data points using the zipper algorithm. At this resolution, it is not
possible to see the difference between these curves. The zipper
algorithm was applied to $6,168$ data points and took 36 seconds. The
image of $24,673$ points on the unit circle took 48 seconds and all of these
points were within $9\cdot 10^{-5}$ of the polygon formed by
connecting the $6,168$ data points. The points on the circle
corresponding to the $6,168$ vertices were mapped to points within
$10^{-10}$ of the verticies.  This error is due to the tolerence set
for Newton's method, round-off error, and the compression/expansion of
harmonic measure.
The image
of $8,160$ verticies in the Carleson grid took 25 seconds to be mapped to
the interior and 25 seconds to the exterior. 

\bigskip
\centerline{\hskip -0.2truein\includegraphics[height=2.3in]{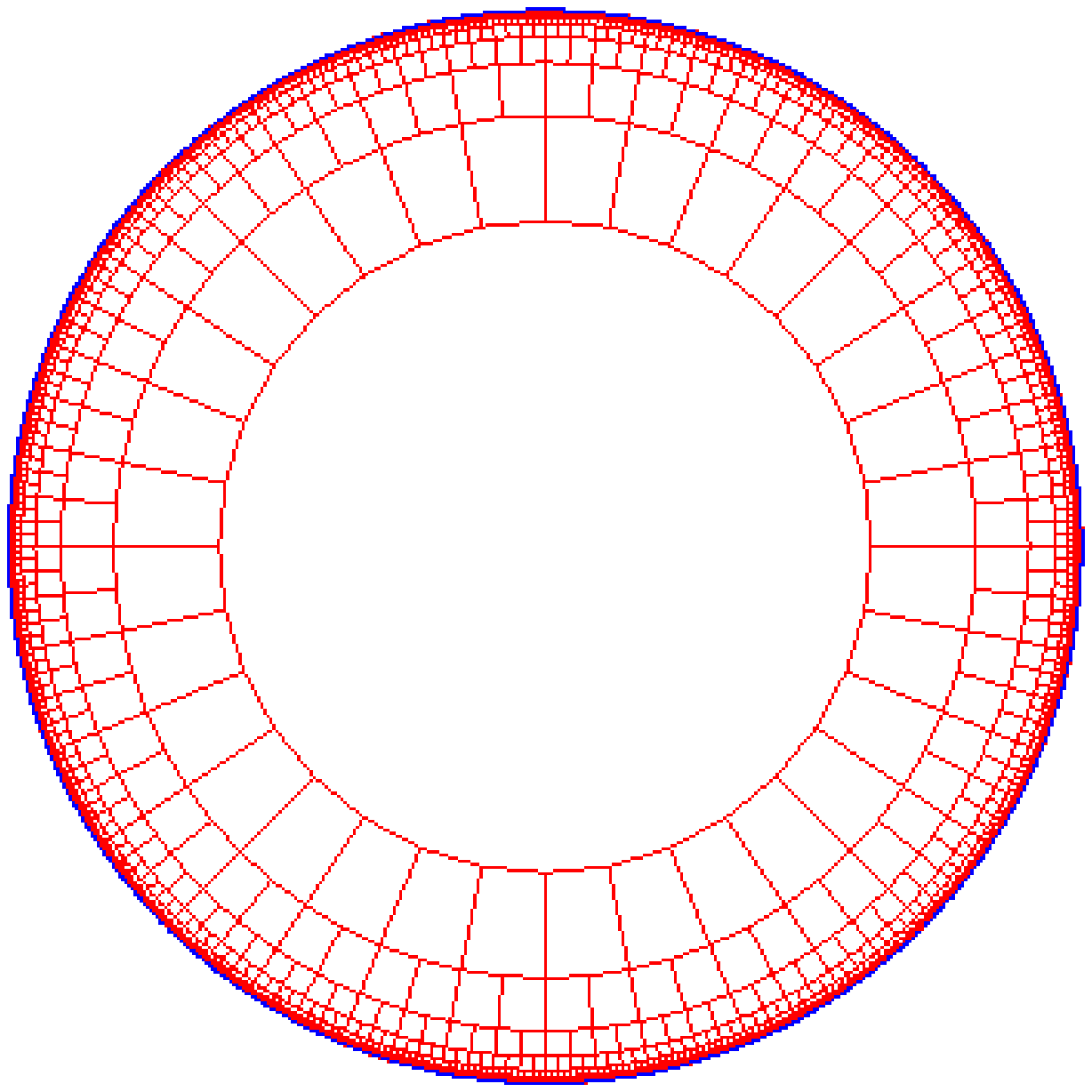} 
\hskip 0.67truein \includegraphics[height=2.4in]{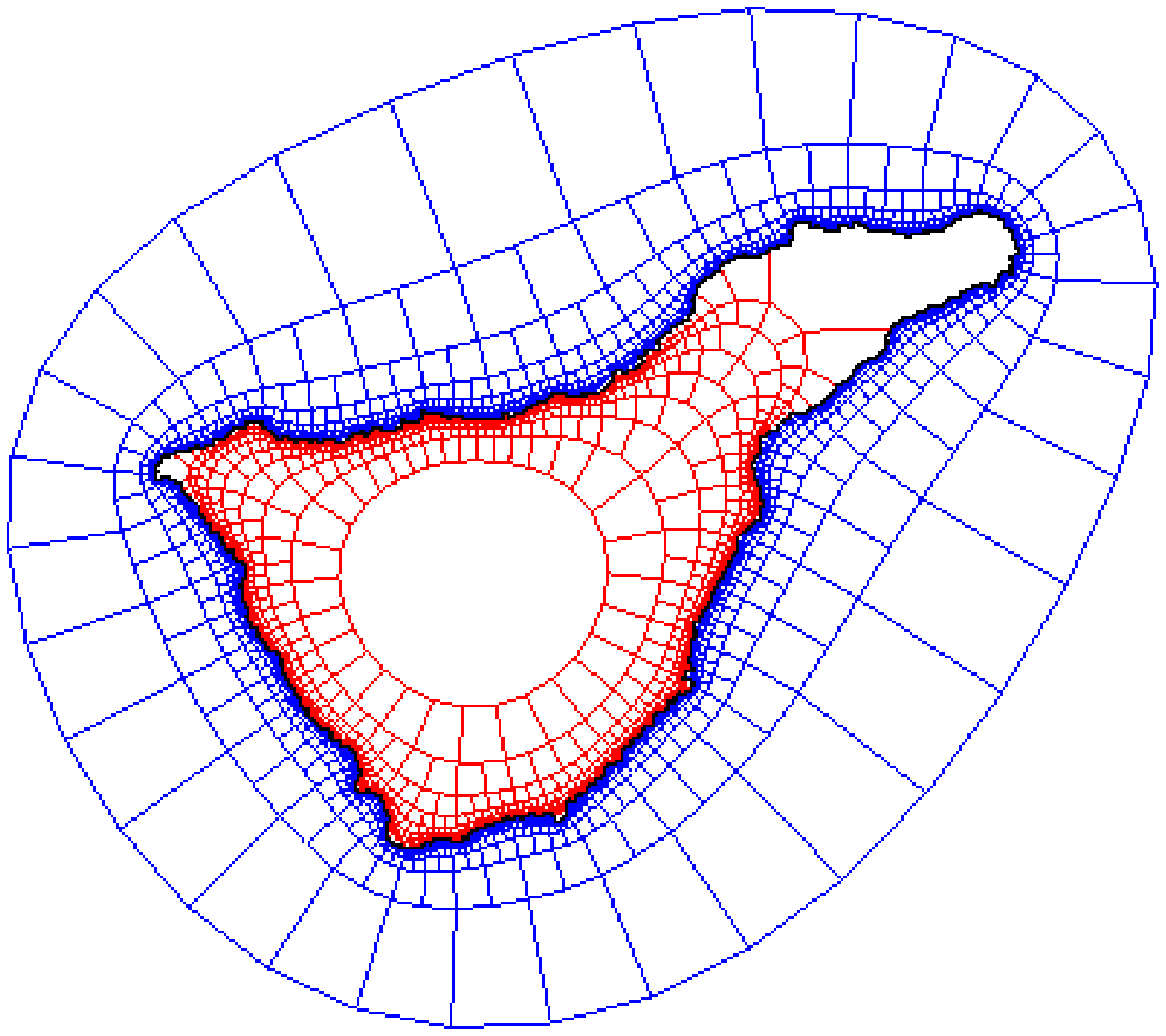}}
\nobreak
\centerline{{\bf Figure 13.} Tenerife.}
\bigskip

Higher resolution images can be obtained from:
\smallskip
http://www.math.washington.edu/$\sim$marshall/preprints/zipper.pdf

The first objection one might have in applying these algorithms with a
large number of data points is that compositions of even very simple
analytic maps can be quite chaotic. Indeed this is the subject of the
field complex dynamics. 
We could redefine the basic maps $f_a$ by composing
with a linear fractional transformation of
the upper half plane so that the composed map is asymptotic to $z$ as
$z\to \infty$. This will not affect the computed curve in these
algorithms since the next basic map begins with a linear fractional
transformation (albeit altered). 
However, if we formulate the basic maps in this way,
then because the maps are nearly linear near $\infty$, the
numerical errors will accumulate only linearly.

Osculation methods also approximate a conformal map by repeated
composition of simple maps.
See Henrici [H] for a discussion of osculation methods and uniform convergence
on compact sets.
The algorithms of the
present article follow the boundary of a given region much 
more closely than, for
instance, the Koebe algorithm and give uniform convergence rather than 
just uniform on compacta.
It is possible to use the techniques of this paper to prove the
geodesic algorithm is an osculation
method for smooth curves, and therefore by the results in [H] 
converge uniformly on
compact subsets. However, prior to this article even a proof that these methods
satisfied the osculation family conditions was not known.

Recently Banjai and Trefethen [BT] adapted multigrid techniques to the
Schwarz-Christoffel algorithm and successfully computed the conformal map to a
region bounded by a polygon with about $10^5$ edges. They used a 12 fold
symmetry in the region
to immediately reduce the parameter problem to size $10^4$. 
Any other conformal mapping
technique can also use symmetry
and obtain a 12 fold reduction in the number of data points required, however  
their work does show at least that Schwarz-Christoffel is possible with $10^4$ vertices, 
though convergence of the algorithm to solve the parameter problem is not always
assured.
The zipper algorithm is competitive in speed and accuracy for such regions. 
The geodesic algorithm is
almost as good, and has the advantage that it is very easy to code and convergence
can be proved. 
It would be interesting to try to prove convergence of the technique used in 
[BT] to find the prevertices, for polygons which are $K$-quasicircles in terms of
$K$.  It would be interesting as well to apply multigrid techniques to 
the zipper algorithm.  

One additional observation worth repeating in this context 
is that the geodesic and zipper algorithms  {\it always} compute a
conformal map of $\H$ to a region bounded by a Jordan curve 
passing through the data points, even if the disc-chain or pacman
conditions are not met.
The image region can be found by evaluating the function at a large number of points
on the real line. By Proposition 2.5 and Corollary 3.9, 
if the data points $\{z_j\}$
satisfy the hypotheses of Theorem 2.2 or Theorem 3.4, then $\varphi$ can be
analytically extended to be a conformal map of the original region
$\Om$ to a region very close to $\D$. 
To do so requires careful consideration of the appropriate
branch of $\sqrt z$ at each stage of the composition.

Theorem 2.2 and Theorem 3.4 and their proofs suggest how to select points on the
boundary of a region to give good accuracy for the mapping functions. Roughly
speaking, points need to be chosen closer together where the region comes close to
folding back on itself. See Figure 12 for example.
Greater accuracy can be obtained by placing more
points on the boundary near the center and fewer on the big lobes. See
also the remarks after the statement of Theorem 6.3 in this regard.
In practice,
the zipper map works well if points are distributed so that 
$$B(z_k, 5|z_{k+1}-z_k|)\cap \b \Om\eqno(7.1)$$
is connected. 
\bigskip
When the boundary of the given region
is not smooth, then one of the processes described in section 2 should
be used to generate the boundary data, if the geodesic algorithm is to
be used. For example, if nothing is
known about the boundary except for a list of data points, then we
preprocess the data by taking data points along the line
segments between the original data points, so that these new points
 correspond to
points of tangency of disjoint circles centered on the line segments,
including circles centered at the original data points.
Note that the original boundary points are not among these new data
points. The geodesic algorithm then finds a conformal map to a region
with the new data points on the boundary. The boundary of the new
region will be close to the polygonal curve through the original data
points, but will not pass through the original data points. This
boundary is ``rounded'' near the original data points. Indeed it is a
smooth curve.

When the boundary of the desired region is less smooth, for example
with ``corners'', then the zipper or slit algorithms should be used. In this
case additional points are placed along the line segments between the
data points, with at least 5 points per edge and satisfying (7.1).
In practice, at least 500 points are chosen on the boundary so that
the image of the circle will be close to the polygonal line through the
data points. Since two data points are pulled down to the real line
with each basic map in the zipper algorithm, 
the original data points should occur at even
numbered indices in the resulting data set (the first data point is
called $z_0$). Then the computed boundary
$\Om_c$ will have corners at each of the original data points, with
angles very close to the angles of the polygon through the original
data points.

\bigskip
A version of the zipper algorithm can be obtained from [MD].
The conformal mapping programs are written in Fortan. Also included is
a graphics program, written in C with X-11 graphics by Mike Stark, 
for the display of the conformal maps. There are also several demo
programs applying the algorithm to problems in elementary fluid flow, 
extremal length and the hyperbolic geometry.
Extensive testing of the geodesic algorithm [MM] and an early version  zipper algorithm 
was done in the 1980's with
Jim Morrow. In particular that experimentation suggested the initial
function $\varphi_0$ in the zipper algorithm which maps the complement of a circular arc
through $z_0$, $z_1$, and $z_2$ onto $\H$.

\bigskip
\bigskip
{\bf Appendix.  J\o rgensen's Theorem.}

\bigskip
Since J\o rgensen's theorem is a key component of the 
proof of the convergence of the geodesic
algorithm, we include a short self-contained proof.  It says that
discs are strictly convex in the hyperbolic geometry of a simply
connected domain $\Om$ (unless $\b\Om$ is contained in the boundary of
the disk).

\bigskip
\proclaim Theorem A.1 (J\o rgensen [J]).
Suppose $\Om$ is a simply connected domain.
If $\Delta$ is
an open disc contained in $\Om$ 
 and if $\gamma$ is a hyperbolic geodesic in $\Om$,
then $\gamma\cap \Delta$ is connected and if non-empty, it is not
tangent to $\b\Delta$ in $\Om$.
\bigskip

\proof See [P, page 91-93]. 
Applying a linear fractional transformation, we may suppose
that the upper half plane $\H \subset \Om$. Suppose $x\in \R$ and
suppose that $f$ is a conformal map of $\D$ onto $\Om$ such that 
$f(0)=x$ and $f'(0)>0$. Then
$$\Im\biggl({{f'(0)}\over {f(z)-x}} -({1\over z} + z)\biggl)$$ 
is a bounded harmonic function on $\D$ which is greater than or equal to
$0$ by the maximum principle. Thus $\Im{{f'(0)}\over{f(z)-x}} \ge 0$ on
$(-1,1)$ and hence $\Im f(z) \le 0$ on the diameter
$(-1,1)$. 
The condition $f'(0)>0$ means that the
geodesic $f\bigl((-1,1)\bigl)$ is tangent to $\R$ at $x$. Thus if
$\gamma$ is a
geodesic which intersects $\H$ and contains the point $x$, then it
cannot be tangent to $\R$ at $x$. Two circles which are orthogonal to
$\R$ can meet in $\H$ in at most one point, and hence hyperbolic geodesics
in simply connected domains (images of orthogonal circles)
meet in at most one point. Thus $\gamma$
cannot reenter $\H$ after leaving it at $x$ because it is separated
from $\R$ by the geodesic $f\bigl((-1,1)\bigl)$. The Theorem follows.
\endproof

In Section 2, we commented that a constructive proof of the Riemann mapping theorem followed from
the proof of Theorem 2.2.  The application of J\o rgensen's theorem in the proof of Theorem 2.2 is
only to domains for which the Riemann map has been explicitly constructed.

\bigskip
\bigskip
{\bf  Bibliography}
\bigskip

\item{[A]} L. Ahlfors, {\it Lectures on quasiconformal mappings}, Van Nostrand (1966).

\item{[BT]} L. Banjai and L. N. Trefethen, {\it A multipole method for Schwarz-Christoffel
mapping of polygons with thousands of sides}, SIAM J. Sci. Comput. 25 (2003)
1042-1065.

\item{[GM]} J. Garnett and D.E. Marshall, {\it Harmonic Measure},
Cambridge Univ. Press (2005).

\item{[H]} P. Henrici, {\it Applied and Computational Complex Analysis}, vol. 3,
J. Wiley  \& Sons(1986).

\item{[J]} J\o rgensen, V., {\it On an inequality for the hyperbolic measure and  its
applications in the theory of functions}, Math. Scand. 4 (1956),
113-124.

\item{[K]} R. K\"uhnau, {\it Numerische Realisierung konformer Abbildungen durch ``Interpolation''}
Z. Angew. Math. Mech. 63 (1983), 631--637 (in German). 

\item{[L]} O. Lehto, {\it Univalent functions and Teichmueller spaces}, Springer (1986).

\item{[M]} A. R. Marchenko, {\it Sur la representation conforme}, C. R. Acad. Sci. USSR vol. 1 (1935), 289--290.

\item{[MD]} D. E. Marshall, {\it Zipper}, 
Fortran programs for numerical computation of conformal
maps, and C programs for X-11 graphics display of the maps.
Sample pictures, Fortran and C code available at URL:

http://www.math.washington.edu/$\sim$marshall/personal.html 

\item{[MM]} D. E. Marshall and J. A. Morrow, {\it Compositions of Slit
Mappings}, unpublished manuscript, 1987.

\item{[P]} Chr. Pommerenke, {\it Boundary behaviour of conformal maps}, Springer (1992).

\item{[SS]}  S. Smale, {\it On the efficiency of algorithms of analysis}.
Bull. Amer. Math. Soc. 13(1985), 87-121.

\item{[SK]} K. Stephenson, {\it Circle packing: a mathematical tale}, Notices
Amer. Math. Soc. 50 (2003), 1376-1388.

\item{[T]} M. Tsuji, {\it Potential theory in modern function theory},
Chelsea (1975).
\item{[W]} S. Warschawski, {\it On the degree of variation in conformal 
mapping of variable regions}, Trans. Amer. Math. Soc. 69 (1950), 335--356.

\end